\documentclass{amsart}
\usepackage{amssymb, amsmath, mathrsfs, verbatim, url, bbm}

\theoremstyle{plain}
\newtheorem{theorem}{Theorem}[section]
\newtheorem{lemma}[theorem]{Lemma}
\newtheorem{proposition}[theorem]{Proposition}
\newtheorem{corollary}[theorem]{Corollary}
\newtheorem*{claim}{Claim}
\newcommand{\MainTheoremName}{Main Theorem}

\theoremstyle{definition}
\newtheorem{definition}[theorem]{Definition}

\theoremstyle{remark}
\newtheorem*{remark}{Remark}

\newtheorem{question}{Question}

\numberwithin{equation}{section}

\DeclareMathOperator{\dom}{dom}
\DeclareMathOperator{\ran}{ran}
\DeclareMathOperator{\cf}{cf}

\newcommand{\nbd}{\nobreakdash}
\newcommand{\homeo}{\cong}
\newcommand{\powset}[1]{\mathcal{P}(#1)}
\newcommand{\card}[1]{\lvert #1\rvert}

\newcommand{\la}{\langle}
\newcommand{\ra}{\rangle}
\newcommand{\ord}{{\mathrm On}}
\newcommand{\carda}{\mathfrak{a}}
\newcommand{\cardb}{\mathfrak{b}}
\newcommand{\cardd}{\mathfrak{d}}
\newcommand{\cardc}{\mathfrak{c}}
\newcommand{\cardg}{\mathfrak{g}}
\newcommand{\cardh}{\mathfrak{h}}
\newcommand{\cardi}{\mathfrak{i}}
\newcommand{\cardp}{\mathfrak{p}}
\newcommand{\cardr}{\mathfrak{r}}
\newcommand{\cards}{\mathfrak{s}}
\newcommand{\cardss}{\cards\cards}
\newcommand{\cardu}{\mathfrak{u}}
\newcommand{\inv}[1]{#1^{-1}}

\newcommand{\restrict}{\upharpoonright}
\newcommand{\MA}[2]{{\mathrm MA}(#1;\, #2)}

\newcommand{\weight}[1]{w(#1)}
\newcommand{\ow}[1]{Nt(#1)}
\newcommand{\owbox}[1]{Nt_\mathrm{box}(#1)}
\newcommand{\opi}[1]{\pi Nt(#1)}
\newcommand{\ochar}[1]{\chi Nt(#1)}
\newcommand{\opichar}[1]{\pi\chi Nt(#1)}

\newcommand{\character}[1]{\chi(#1)}

\newcommand{\stonecech}{\beta}

\newcommand{\Fn}[2]{\mathrm{Fn}(#1,\,#2)}
\newcommand{\wma}{we may assume}
\newcommand{\Wma}{We may assume}
\newcommand{\op}{\mathrm{op}}

\newcommand{\omegaop}{$\omega^\op$\nbd-like}
\newcommand{\kappaop}{$\kappa^\op$\nbd-like}
\newcommand{\cardop}[1]{$#1^\op$\nbd-like}

\newcommand{\forces}{\Vdash}
\newcommand{\supp}[1]{\mathrm{supp}(#1)}

\newcommand{\additivity}[1]{\mathrm{add}(#1)}

\newcommand{\pol}[2]{{#1\restrict #2}}
\newcommand{\poleq}[2]{#1\restrict_{\leq} #2}
\newcommand{\fod}[2]{{#1\restrict #2}}
\newcommand{\fol}[2]{{#1\restrict #2}}
\newcommand{\foleq}[2]{#1\restrict_{\leq} #2}
\newcommand{\gfd}[2]{{#1\restrict #2}}
\newcommand{\gfl}[2]{{#1\restrict #2}}
\newcommand{\gfleq}[2]{#1\restrict_{\leq} #2}
\newcommand{\gfp}[2]{{#1}_{#2}}
\newcommand{\fnl}[2]{\pol{#1}{#2}}
\newcommand{\fnleq}[2]{\poleq{#1}{#2}}

\newcommand{\mcA}{\mathcal{ A}}
\newcommand{\mcB}{\mathcal{ B}}
\newcommand{\mcC}{\mathcal{ C}}
\newcommand{\mcD}{\mathcal{ D}}
\newcommand{\mcE}{\mathcal{ E}}
\newcommand{\mcF}{\mathcal{ F}}

\newcommand{\mcI}{\mathcal{ I}}
\newcommand{\mcJ}{\mathcal{ J}}

\newcommand{\mcU}{\mathcal{ U}}
\newcommand{\mcV}{\mathcal{ V}}

\newcommand{\mbP}{\mathbb{ P}}
\newcommand{\mbQ}{\mathbb{ Q}}
\newcommand{\mbR}{\mathbb{ R}}

\newcommand{\meager}{\mcB}

\begin{document}

\title{Splitting families and the Noetherian type of $\stonecech\omega\setminus\omega$}

\author{David Milovich}

\date{March 26, 2008}

\address{University of Wisconsin-Madison Mathematics Dept.}
\email{milovich@math.wisc.edu}
\thanks{Support provided by an NSF graduate fellowship.}

\begin{abstract}
Extending some results of Malykhin, we prove several independence results about base properties of $\stonecech\omega\setminus\omega$ and its powers, especially the Noetherian type $\ow{\stonecech\omega\setminus\omega}$, the least $\kappa$ for which $\stonecech\omega\setminus\omega$ has a base that is $\kappa$\nbd-like with respect to containment.  For example, $\ow{\stonecech\omega\setminus\omega}$ is at least $\cards$, but can consistently be $\omega_1$, $\cardc$, $\cardc^+$, or strictly between $\omega_1$ and $\cardc$.  $\ow{\stonecech\omega\setminus\omega}$ is also consistently less than the additivity of the meager ideal.  $\ow{\stonecech\omega\setminus\omega}$ is closely related to the existence of special kinds of splitting families.
\end{abstract}

\maketitle

\section{Introduction}

\begin{definition}
Given a cardinal $\kappa$, define a poset to be $\kappa$\nbd-\emph{like} ($\kappa^{\op}$\nbd-\emph{like}) if no element is above (below) $\kappa$\nbd-many elements.  Define a poset to be \emph{almost} $\kappa^{\op}$\nbd-\emph{like} if it has a \kappaop\ dense subset.
\end{definition}

In the context of families of subsets of a topological space, we will always implicitly order by inclusion.  We are particularly interested in \kappaop\ bases, $\pi$\nbd-bases, local bases, and local $\pi$\nbd-bases of the space $\omega^*$ of nonprincipal ultrafilters on $\omega$.  Recall that a local base (local $\pi$\nbd-base) at a point in a space is a family of open neighborhoods of that point (family of nonempty open subsets) such that every neighborhood of the point contains an element of the family; a base ($\pi$\nbd-base) of a space is family of open sets that contains local bases (local $\pi$\nbd-bases) at every point.  See Engelking~\cite{engelking} for the more background on bases and their cousins.  Also recall the following basic cardinal functions.  For more about these functions, see Juh\'asz~\cite{juhasz}.

\begin{definition}
Given a space $X$, let the \emph{weight} of $X$, or $\weight{X}$, be the least $\kappa\geq\omega$ such that $X$ has a base of size at most $\kappa$.  Given $p\in X$, let the \emph{character} of $p$, or $\character{p,X}$, be the least $\kappa\geq\omega$ such that there is a local base at $p$ of size at most $\kappa$.  Let the character of $X$, or $\character{X}$, be the supremum of the characters of its points.  Analogously define \emph{$\pi$\nbd-weight} and \emph{local $\pi$\nbd-character}, respectively denoting them using $\pi$ and $\pi\chi$.
\end{definition}

Now consider the following order\nbd-theoretic parallels.

\begin{definition}
Given a space $X$, let the \emph{Noetherian type} of $X$, or $\ow{X}$, be the least $\kappa\geq\omega$ such that $X$ has a base that is \kappaop.  Given $p\in X$, let the \emph{local Noetherian type} of $p$, or $\ochar{p,X}$, be the least $\kappa\geq\omega$ such that there is a \kappaop\ local base at $p$.  Let the local Noetherian type of $X$, or $\ochar{X}$, be the supremum of the local Noetherian types of its points.  Analogously define \emph{Noetherian $\pi$\nbd-type} and \emph{local Noetherian $\pi$\nbd-type}, respectively denoting them using $\pi Nt$ and $\pi\chi Nt$.
\end{definition}

Noetherian type and Noetherian $\pi$\nbd-type were introduced by Peregudov~\cite{peregudov97}.  Let $\omega^*$ denote the space of nonprincipal ultrafilters on $\omega$.  Malykhin~\cite{malykhin} proved that MA implies $\opi{\omega^*}=\cardc$ and CH implies $\ow{\omega^*}=\cardc$.  We extend these results by investigating $\ow{\omega^*}$, $\opi{\omega^*}$, $\ochar{\omega^*}$, and $\opichar{\omega^*}$ as cardinal characteristics of the continuum.  For background on such cardinals, see Blass~\cite{blass}.  We also examine the sequence $\la\ow{(\omega^*)^{1+\alpha}}\ra_{\alpha\in\ord}$.

\begin{definition}  Let $\cardb$ denote the minimum of $\card{\mcF}$ where $\mcF$ ranges over the subsets of $\omega^\omega$ that have no upper bound in $\omega^\omega$ with respect to eventual domination.
\end{definition}

\begin{definition}
A \emph{tree $\pi$\nbd-base} of a space $X$ is a $\pi$\nbd-base that is a tree when ordered by containment.  Let $\cardh$ be the minimum of the set of heights of tree $\pi$\nbd-bases of $\omega^*$.
\end{definition}

Balcar, Pelant, and Simon~\cite{balcarcardh} proved that tree $\pi$\nbd-bases of $\omega^*$ exist, and that $\cardh\leq\min\{\cardb,\cf\cardc\}$.  They also proved that the above definition of $\cardh$ is equivalent to the more common definition of $\cardh$ as the distributivity number of $[\omega]^{\omega}$ ordered by $\subseteq^*$.

\begin{definition}
Given $x,y\in[\omega]^\omega$, we say that $x$ \emph{splits} $y$ if $\card{y\cap x}=\card{y\setminus x}=\omega$.  Let $\cardr$ be the minimum value of $\card{A}$ where $A$ ranges over the subsets of $[\omega]^\omega$ such that no $x\in[\omega]^\omega$ splits every $y\in A$.  Let $\cards$ be the minimum value of $\card{A}$ where $A$ ranges over the subsets of $[\omega]^\omega$ such that every $x\in[\omega]^\omega$ is split by some $y\in A$.
\end{definition}

It is known that $\cardb\leq\cardr$ and $\cardh\leq\cards$. (See Theorems 3.8 and 6.9 of \cite{blass}.)

Clearly, $\ow{\omega^*}\leq\weight{\omega^*}^+=\cardc^+$.  We will show that also $\opichar{\omega^*}=\omega$ and $\opi{\omega^*}=\cardh$ and $\cards\leq\ow{\omega^*}$.  Furthermore, $\ow{\omega^*}$ can consistently be $\cardc$, $\cardc^+$, or any regular $\kappa$ satisfying $2^{<\kappa}=\cardc$.  Also, $\ow{\omega^*}=\omega_1$ is relatively consistent with any values of $\cardb$ and $\cardc$.  The relations $\omega_1<\cardb=\cards=\ow{\omega^*}<\cardc$ and $\omega_1=\cardb=\cards<\ow{\omega^*}<\cardc$ are also each consistent.  We also prove some relations between $\cardr$ and $\ow{\omega^*}$, as well as some consistency results about the local Noetherian type of points in $\omega^*$.

\section{Basic results}

The following proposition is essentially due to Peregudov (see Lemma 1 of \cite{peregudov97}).

\begin{proposition}\label{PROpicharcfchar}
Suppose a point $p$ in a space $X$ satisfies $\pi\character{p,X}<\cf\kappa\leq\kappa\leq\character{p,X}$.  Then $\ow{X}>\kappa$.
\end{proposition}
\begin{proof}
Let $\mcA$ be a base of $X$.  Let $\mcU_0$ and $\mcV_0$ be, respectively, a local $\pi$\nbd-base at $p$ of size at most $\pi\character{p,X}$ and a local base at $p$ of size $\character{p,X}$.  For each element of $\mcU_0$, choose a subset in $\mcA$, thereby producing local $\pi$\nbd-base $\mcU$ at $p$ that is a subset of $\mcA$ of size at most $\pi\character{p,X}$.  Similarly, for each element of $\mcV_0$, choose a smaller neighborhood of $p$ in $\mcA$, thereby producing a local base $\mcV$ at $p$ that is a subset of $\mcA$ of size $\character{p,X}$.  Every element of $\mcV$ contains an element of $\mcU$.  Hence, some element of $\mcU$ is contained in $\kappa$\nbd-many elements of $\mcV$; hence, $\mcA$ is not \kappaop.
\end{proof}

\begin{definition}
For all $x\in[\omega]^\omega$, set $x^*=\{p\in\omega^*:p\in x\}$.
\end{definition}

\begin{theorem}\label{THMowomegastarcplus}
It is relatively consistent with any value of $\cardc$ satisfying $\cf\cardc>\omega_1$ that $\ow{\omega^*}=\cardc^+$.
\end{theorem}
\begin{proof}
\Wma\ $\cf\cardc>\omega_1$.  By Exercise A10 on p. 289 of Kunen~\cite{kunen}, there is a ccc generic extension $V[G]$ such that $\check{\cardc}=\cardc^{V[G]}$ and, in $V[G]$, there exists $p\in\omega^*$ such that $\character{p,\omega^*}=\omega_1$.  Henceforth work in $V[G]$.  Let $\varphi$ be a bijection from $\omega^2$ to $\omega$.  Define $\psi\colon\omega^*\to\omega^*$ by 
\begin{equation*}
x\mapsto\{E\subseteq\omega:\{m<\omega:\{n<\omega:\varphi(m,n)\in E\}\in p\}\in x\}.
\end{equation*}
Since $\pi\character{p,\omega^*}\leq\character{p,\omega^*}=\omega_1$, there exists $\la E_\alpha\ra_{\alpha<\omega_1}\in([\omega]^\omega)^{\omega_1}$ such that every neighborhood of $p$ contains $E_\alpha^*$ for some $\alpha<\omega_1$.  Hence, for all $x\in\omega^*$, every neighborhood of $\psi(x)$ contains $(\varphi``(\{m\}\times E_\alpha))^*$ for some $m<\omega$ and $\alpha<\omega_1$; whence, $\pi\character{\psi(x),\omega^*}=\omega_1$.  Since $\psi$ is easily verified to be a topological embedding, $\character{x,\omega^*}\leq\character{\psi(x),\omega^*}$ for all $x\in\omega^*$.  By a result of Pospi\v sil~\cite{pospisil}, there exists $q\in\omega^*$ such that $\character{q,\omega^*}=\cardc$.  Hence, $\pi\character{\psi(q),\omega^*}=\omega_1$ and $\character{\psi(q),\omega^*}=\cardc$.  By Proposition~\ref{PROpicharcfchar}, $\ow{\omega^*}>\character{\psi(q),\omega^*}=\cardc$.
\end{proof}

\begin{definition}
Given $n<\omega$, let $\cardss_n$ ($\cardss_\omega$) denote the least cardinal $\kappa$ for which there exists a sequence $\la f_\alpha\ra_{\alpha<\cardc}$ of functions on $\omega$ each with range contained in $n$ (each with finite range) such that for all $I\in[\cardc]^\kappa$ and $x\in[\omega]^\omega$ there exists $\alpha\in I$ such that $f_\alpha$ is not eventually constant on $x$.  (The notation $\cardss$ was chosen with the phrase ``supersplitting number'' in mind.)  Note that if such an $\la f_\alpha\ra_{\alpha<\cardc}$ does not exist for any $\kappa\leq\cardc$, then $\cardss_n$ ($\cardss_\omega$) is by definition equal to $\cardc^+$.
\end{definition}

Clearly $\cardss_n\geq\cardss_{n+1}\geq\cardss_\omega$ for all $n<\omega$.  Moreover, since $\cf\cardc>\omega$, we have $\cardss_\omega=\cardss_n$ for some $n<\omega$.  However, for any particular $n\in\omega\setminus 2$, it is not clear whether ZFC proves $\cardss_\omega=\cardss_n$.

\begin{definition}
Given $\lambda\geq\kappa\geq\omega$ and a space $X$, a $\la\lambda,\kappa\ra$\nbd-\emph{splitter} of $X$ is a sequence $\la\mcF_\alpha\ra_{\alpha<\lambda}$ of finite open covers of $X$ such that, for all $I\in[\lambda]^\kappa$ and $\la U_\alpha\ra_{\alpha\in I}\in\prod_{\alpha\in I}\mcF_\alpha$, the interior of $\bigcap_{\alpha\in I}U_\alpha$ is empty. 
\end{definition}

\begin{lemma}\label{LEMowsplitterupperbound}
Suppose $X$ is a compact space with a base $\mcA$ of size at most $\weight{X}$ such that $U\cap V\in\mcA\cup\{\emptyset\}$ for all $U,V\in\mcA$.  If $\kappa\leq\weight{X}$ and $X$ has a $\la\weight{X},\kappa\ra$\nbd-splitter, then $\mcA$ contains a \kappaop\ base of $X$.  Hence, $\ow{\omega^*}\leq\cardss_\omega$.
\end{lemma}
\begin{proof}
Set $\lambda=\weight{X}$ and let $\la\mcF_\alpha\ra_{\alpha<\lambda}$ be a $\la\lambda,\kappa\ra$\nbd-splitter of $X$.  For each $\alpha<\lambda$, the cover $\mcF_\alpha$ is refined by a finite subcover of $\mcA$; hence, \wma\ $\mcF_\alpha\subseteq\mcA$.  Let $\mcA=\{U_\alpha:\alpha<\lambda\}$.  For each $\alpha<\lambda$, set $\mcB_\alpha=\{U_\alpha\cap V: V\in\mcF_\alpha\}$.  Set $\mcB=\bigcup_{\alpha<\lambda}\mcB_\alpha\setminus\{\emptyset\}$.  Then $\mcB$ is easily seen to be a base of $X$ and a \kappaop\ subset of $\mcA$.
\end{proof}

\begin{lemma}\label{LEMowlowerboundcharweight}
Let $X$ be a compact space without isolated points and let $\omega\leq\kappa\leq\lambda\leq\min_{p\in X}\character{p,X}$.  If $X$ has no $\la\lambda,\kappa\ra$\nbd-splitter, then $\ow{X}>\kappa$.
\end{lemma}
\begin{proof}
Let $\mcA$ be a base of $X$.  Construct a sequence $\la\mcF_\alpha\ra_{\alpha<\lambda}$ of finite subcovers of $\mcA$ as follows.  Suppose we have $\alpha<\lambda$ and $\la\mcF_\beta\ra_{\beta<\alpha}$.  For each $p\in X$, choose $V_p\in\mcA$ such that $p\in V_p\not\in\bigcup_{\beta<\alpha}\mcF_\beta$.  Let $\mcF_\alpha$ be a finite subcover of $\{V_p:p\in X\}$.  Then $\mcF_\alpha\cap\mcF_\beta=\emptyset$ for all $\alpha<\beta<\lambda$.  Suppose $X$ has no $\la\lambda,\kappa\ra$\nbd-splitter.  Then choose $I\in[\lambda]^\kappa$ and $\la U_\alpha\ra_{\alpha\in I}\in\prod_{\alpha\in I}\mcF_\alpha$ such that $\bigcap_{\alpha\in I}U_\alpha$ has nonempty interior.  Then there exists $W\in\mcA$ such that $W\subseteq\bigcap_{\alpha\in I}U_\alpha$.  Thus, $\mcA$ is not \kappaop.
\end{proof}

\begin{definition}
Let $\cardu$ denote the minimum of the set of characters of points in $\omega^*$.  Let $\pi\cardu$ denote the minimum of the set of $\pi$\nbd-characters of points in $\omega^*$.  
\end{definition}

By a theorem of Balcar and Simon~\cite{balcarcardr}, $\pi\cardu=\cardr$.

\begin{theorem}\label{THMcarducardc}
Suppose $\cardu=\cardc$.  Then $\ow{\omega^*}=\cardss_\omega$.
\end{theorem}
\begin{proof}
By Lemma~\ref{LEMowsplitterupperbound}, $\ow{\omega^*}\leq\cardss_\omega$.  Suppose $\kappa\leq\cardc$.  Since every finite open cover of $\omega^*$ is refined by a finite, pairwise disjoint, clopen cover, $\omega^*$ has a $\la\cardc,\kappa\ra$\nbd-splitter if and only if $\cardss_\omega\leq\kappa$.  Hence, $\ow{\omega^*}\geq\cardss_\omega$ by Lemma~\ref{LEMowlowerboundcharweight}.
\end{proof}

\begin{lemma}\label{LEMcrss2}
Suppose $\cardr=\cardc$.  Then $\cardss_2\leq\cardc$.
\end{lemma}
\begin{proof}
Let $\la x_\alpha\ra_{\alpha<\cardc}$ enumerate $[\omega]^\omega$.  Construct $\la y_\alpha\ra_{\alpha<\cardc}\in([\omega]^\omega)^\cardc$ as follows.  Given $\alpha<\cardc$ and $\la y_\beta\ra_{\beta<\alpha}$, choose $y_\alpha$ such that $y_\alpha$ splits every element of $\{x_\alpha\}\cup\{y_\beta:\beta<\alpha\}$.  Suppose $I\in[\cardc]^\cardc$ and $\alpha<\cardc$.  Then $x_\alpha$ is split by $y_\beta$ for all $\beta\in I\setminus\alpha$.  Thus, $\la\{y_\alpha,\omega\setminus y_\alpha\}\ra_{\alpha<\cardc}$ witnesses $\cardss_2\leq\cardc$.
\end{proof}

\begin{theorem}\label{THMowcardr}  The cardinals $\cardr$ and $\ow{\omega^*}$ are related as follows.
\begin{enumerate}
\item If $\cardr=\cardc$, then $\ow{\omega^*}=\cardss_\omega\leq\cardc$.
\item If $\cardr<\cardc$, then $\ow{\omega^*}\geq\cardc$.  
\item If $\cardr<\cf\cardc$, then $\ow{\omega^*}=\cardc^+$.
\end{enumerate}
\end{theorem}
\begin{proof}
Statement (1) follows from Lemma~\ref{LEMcrss2}, Theorem~\ref{THMcarducardc}, and $\pi\cardu=\cardr$.  The proof of Theorem~\ref{THMowomegastarcplus} shows how to construct $p\in\omega^*$ such that $\pi\character{p,\omega^*}=\pi\cardu=\cardr$ and $\character{p,\omega^*}=\cardc$.  Hence, (2) and (3) follow from Proposition~\ref{PROpicharcfchar}.
\end{proof}

\begin{definition}
A subset $A$ of $[\omega]^\omega$ has the \emph{strong finite intersection property} (SFIP) if the intersection of every finite subset of $A$ is infinite.  Given $A\subseteq[\omega]^\omega$ with the SFIP, define the \emph{Booth forcing for} $A$ to be $[\omega]^{<\omega}\times[A]^{<\omega}$ ordered by $\la\sigma_0,F_0\ra\leq\la\sigma_1,F_1\ra$ if and only if $F_0\supseteq F_1$ and  $\sigma_1\subseteq\sigma_0\subseteq\sigma_1\cup\bigcap F_1$.  Define a \emph{generic pseudointersection} of $A$ to be $\bigcup_{\la\sigma,F\ra\in G}\sigma$ where $G$ is a generic filter of $[\omega]^{<\omega}\times[A]^{<\omega}$.
\end{definition}

\begin{theorem}
For all cardinals $\kappa$ satisfying $\kappa>\cf\kappa>\omega$, it is consistent that $\cardr=\cardu=\cf\kappa$ and $\ow{\omega^*}=\cardss_2=\cardc=\kappa$.
\end{theorem}
\begin{proof}
Assuming GCH in the ground model, construct a finite support iteration $\la\mbP_\alpha\ra_{\alpha\leq\kappa}$ as follows.  First choose some $U_0\in\omega^*$.  Then suppose we have $\alpha<\kappa$ and  $\mbP_\alpha$ and $\forces_\alpha U_\alpha\in\omega^*$.  Let $\mbP_{\alpha+1}\cong\mbP_\alpha*\mbQ_\alpha$ where $\mbQ_\alpha$ is a $\mbP_\alpha$\nbd-name for the Booth forcing for $U_\alpha$.  Let $x_\alpha$ be a $\mbP_{\alpha+1}$\nbd-name for a generic pseudointersection of $U_\alpha$ added by $\mbQ_\alpha$; let $U_{\alpha+1}$ be a  $\mbP_{\alpha+1}$\nbd-name for an element of $\omega^*$ containing $U_\alpha\cup\{x_\alpha\}$.  For limit $\alpha<\kappa$, let $U_\alpha=\bigcup_{\beta<\alpha}U_\beta$.

Let $\la\eta_\alpha\ra_{\alpha<\cf\kappa}$ be an increasing sequence of ordinals with supremum $\kappa$.  Then $\{x_{\eta_\alpha}:\alpha<\cf\kappa\}$ is forced to generate an ultrafilter in $V^{\mbP_\kappa}$.  Hence, $\forces_\kappa\cardr\leq\cardu\leq\cf\kappa<\kappa=\cardc$.  Therefore, by Lemma~\ref{LEMowsplitterupperbound} and Theorem~\ref{THMowcardr}, it suffices to show that $\forces_\kappa\cardss_2\leq\kappa$.  Every nontrivial finite support iteration of infinite length adds a Cohen real.  Hence, we may choose for each $\alpha<\kappa$ a $\mbP_{\omega(\alpha+1)}$\nbd-name $y_\alpha$ for an element of $[\omega]^\omega$ that is Cohen over $V^{\mbP_{\omega\alpha}}$.  Then every name $S$ for the range of a cofinal subsequence of $\la y_\alpha\ra_{\alpha<\kappa}$ is such that
\begin{equation*}
\forces_\kappa\forall z\in[\omega]^\omega\ \,\exists w\in S\ \ w\text{ splits }z.
\end{equation*}
Hence, $\la y_\alpha\ra_{\alpha<\kappa}$ witnesses that  $\forces_\kappa\cardss_2\leq\kappa$.
\end{proof}

\begin{theorem}\label{THMowcards}
$\ow{\omega^*}\geq\cards$.
\end{theorem}
\begin{proof}
Suppose $\ow{\omega^*}=\kappa<\cards$.  Since $\ow{\omega^*}<\cardc$, we have $\cardr=\cardc$ by Theorem~\ref{THMowcardr}.  Hence, $\cardu=\cardc$.  By Theorem~\ref{THMcarducardc}, it suffices to show that $\cardss_\omega>\kappa$.  Suppose $\la f_\alpha\ra_{\alpha<\cardc}$ is a sequence of functions on $\omega$ with finite range and $I\in[\cardc]^\kappa$.  Since $\kappa<\cards$, there exists $x\in[\omega]^\omega$ such that $f_\alpha$ is eventually constant on $x$ for all $\alpha\in I$.  Thus, $\cardss_\omega>\kappa$.
\end{proof}

\begin{lemma}\label{LEMmutuallydense}
Let $\kappa$ be a cardinal and let $P$ and $Q$ be mutually dense subsets of a common poset.  Then $P$ is almost \kappaop\ if and only if $Q$ is.
\end{lemma}
\begin{proof}
Suppose $D$ is a \kappaop\ dense subset of $P$.  Then it suffices to construct a \kappaop\ dense subset of $Q$.  Define a partial map $f$ from $\card{D}^+$ to $Q$ as follows.  Set $f_0=\emptyset$.  Suppose $\alpha<\card{D}^+$ and we have constructed a partial map $f_\alpha$ from $\alpha$ to $Q$.  Set $E=\{d\in D: d\not\geq q\text{ for all }q\in\ran f_\alpha\}$.  If $E=\emptyset$, then set $f_{\alpha+1}=f_\alpha$.  Otherwise, choose $q\in Q$ such that $q\leq e$ for some $e\in E$ and let $f_{\alpha+1}$ be the smallest function extending $f_\alpha$ such that $f_{\alpha+1}(\alpha)=q$.  For limit ordinals $\gamma\leq\card{D}^+$, set $f_\gamma=\bigcup_{\alpha<\gamma}f_\alpha$.  Set $f=f_{\card{D}^+}$.

Let us show that $\ran f$ is a \kappaop.  Suppose otherwise.  Then there exists $q\in\ran f$ and an increasing sequence $\la\xi_\alpha\ra_{\alpha<\kappa}$ in $\dom f$ such that $q\leq f(\xi_\alpha)$ for all $\alpha<\kappa$.  By the way we constructed $f$, there exists $\la d_\alpha\ra_{\alpha<\kappa}\in D^\kappa$ such that $f(\xi_\beta)\leq d_\beta\not=d_\alpha$ for all $\alpha<\beta<\kappa$.  Choose $p\in P$ such that $p\leq q$.  Then choose $d\in D$ such that $d\leq p$.  Then $d\leq d_\beta\not=d_\alpha$ for all $\alpha<\beta<\kappa$, which contradicts that $D$ is \kappaop.  Therefore, $\ran f$ is \kappaop.

Finally, let us show that $\ran f$ is a dense subset of $Q$.  Suppose $q\in Q$.  Choose $p\in P$ such that $p\leq q$.  Then choose $d\in D$ such that $d\leq p$.  By the way we constructed $f$, there exists $r\in\ran f$ such that $r\leq d$; hence, $r\leq q$.
\end{proof}

\begin{theorem}\label{THMopicardh}
$\opi{\omega^*}=\cardh$.
\end{theorem}
\begin{proof}
First, we show that $\opi{\omega^*}\leq\cardh$.  Let $\mcA$ be a tree $\pi$\nbd-base of $\omega^*$ such that $\mcA$ has height $\cardh$ with respect to containment.  Then $\mcA$ is clearly \cardop{\cardh}.  To show that $\cardh\leq\opi{\omega^*}$, let $\mcA$ be as above and let $\mcB$ be a \cardop{\opi{\omega^*}}\ $\pi$\nbd-base of $\omega^*$.  Then $\mcA$ and $\mcB$ are mutually dense; hence, by Lemma~\ref{LEMmutuallydense}, $\mcA$ contains a \cardop{\opi{\omega^*}}\  $\pi$\nbd-base $\mcC$ of $\omega^*$.  Since $\mcC$ is also a tree $\pi$\nbd-base, it has height at most $\opi{\omega^*}$.  Hence, $\cardh\leq\opi{\omega^*}$.
\end{proof}

\begin{corollary}\label{CORowomegastarc}
If $\cardh=\cardc$, then $\opi{\omega^*}=\ow{\omega^*}=\cardss_2=\cardc$.
\end{corollary}
\begin{proof}
Suppose $\cardh=\cardc$.  Then $\cardr=\cardc$ because $\cardh\leq\cardb\leq\cardr\leq\cardc$.  Hence, by Theorem~\ref{THMopicardh}, Theorem~\ref{THMowcardr}, and Lemma~\ref{LEMcrss2}, $\cardc\leq\opi{\omega^*}\leq\ow{\omega^*}=\cardss_\omega\leq\cardss_2\leq\cardc$.
\end{proof}

\section{Models of $\ow{\omega^*}=\omega_1$}

Adding $\cardc$\nbd-many Cohen reals collapses $\cardss_2$ to $\omega_1$.  By Lemma~\ref{LEMowsplitterupperbound}, it therefore also collapses $\ow{\omega^*}$ to $\omega_1$.  The same result holds for random reals and Hechler reals.

\begin{theorem}\label{THMowrandomcohen}
Suppose $\kappa^\omega=\kappa$ and $\mbP=\mcB(2^\kappa)/\mcI$ where $\mcB(2^\kappa)$ is the Borel alegebra of the product space $2^\kappa$ and $\mcI$ is either the meager ideal or the null ideal (with respect to the product measure).  (In other words, $\mbP$ adds $\kappa$\nbd-many Cohen reals or $\kappa$\nbd-many random reals in the usual way.)  Then $\mathbbm{1}_\mbP\forces\omega_1=\cardss_2$.
\end{theorem}
\begin{proof}
Working in the generic extension $V[G]$, we have $\kappa=\cardc$ and a sequence $\la x_\alpha\ra_{\alpha<\kappa}$ in $[\omega]^\omega$ such that $V[G]=V[\la x_\alpha\ra_{\alpha<\kappa}]$ and, if $E\in\powset{\kappa}\cap V$ and $\alpha\in\kappa\setminus E$, then $x_\alpha$ is Cohen or random over $V[\la x_\beta\ra_{\beta\in E}]$.  (See \cite{kunenrandom} for a proof.)  Suppose $I\in[\kappa]^{\omega_1}$ and $y\in[\omega]^\omega$.  Then $y\in V[\la x_\alpha\ra_{\alpha\in J}]$ for some $J\in[\kappa]^\omega\cap V$; hence, $x_\alpha$ splits $y$ for all $\alpha\in I\setminus J$.  Thus, $\la\{x_\alpha,\omega\setminus x_\alpha\}\ra_{\alpha<\kappa}$ witnesses $\cardss_2=\omega_1$.
\end{proof}

\begin{definition}
Let $\cardd$ denote the minimum of the cardinalities of subsets of $\omega^\omega$ that are cofinal with respect to eventual domination.
\end{definition}

\begin{corollary}
Every transitive model of ZFC has a ccc forcing extension that preserves $\cardb$, $\cardd$, and $\cardc$, and collapses $\cardss_2$ to $\omega_1$.
\end{corollary}
\begin{proof}
Add $\cardc$\nbd-many random reals to the ground model.  Then every element of $\omega^\omega$ in the extension is eventually dominated by an element of $\omega^\omega$ in the ground model; hence, $\cardb$, $\cardd$, and $\cardc$ are preserved by this forcing, while $\cardss_2$ becomes $\omega_1$.
\end{proof}

\begin{definition}
We say that a transfinite sequence $\la x_\alpha\ra_{\alpha<\eta}$ of subsets of $\omega$ is \emph{eventually splitting} if for all $y\in[\omega]^\omega$ there exists $\alpha<\eta$ such that for all $\beta\in\eta\setminus\alpha$ the set $x_\beta$ splits $y$.
\end{definition}

\begin{theorem}\label{THMowhechler}
Let $\kappa=\kappa^\omega$.  Then $\cardss_2=\omega_1$ is forced by the $\kappa$\nbd-long finite support iteration of Hechler forcing.
\end{theorem}
\begin{proof}
Let $\mbP$ be the $\kappa$\nbd-long finite support iteration of Hechler forcing.  Let $G$ be a generic filter of $\mbP$.  For each $\alpha<\kappa$, let $g_\alpha$ be the generic dominating function added at stage $\alpha$; set $x_\alpha=\{n<\omega:g_\alpha(n)\text{ is even}\}$.  Suppose $p\in G$ and $I$ and $y$ are names such that $p$ forces $I\in[\kappa]^{\omega_1}$ and $y\in[\omega]^\omega$.  Choose $q\in G$ and a name $h$ such that $q\leq p$ and $q$ forces $h$ to be an increasing map from $\omega_1$ to $I$.  For each $\alpha<\omega_1$, set $E_\alpha=\{\beta<\kappa:q\not\forces h(\alpha)\not=\check{\beta}\}$; let $k_\alpha$ be a surjection from $\omega$ to $E_\alpha$.  Let $q\geq r\in G$ and $n<\omega$ and $\gamma\leq\kappa$ and $J$ be a name such that $r$ forces $J\in[\omega_1]^{\omega_1}$ and $\sup\ran h=\check{\gamma}$ and $h(\alpha)=k_\alpha(n)^{\check{}}$ for all $\alpha\in J$.  Set $F=\{k_\alpha(n):\alpha<\omega_1\}\cap\gamma$; let $j$ be the order isomorphism from some ordinal $\eta$ to $F$.  Then $\cf\eta=\cf\gamma=\omega_1$.  For all $\alpha<\kappa$, the set $x_\alpha$ is Cohen over $V[\la g_\beta\ra_{\beta<\alpha}]$; hence, $\la x_{j(\alpha)}\ra_{\alpha<\eta}$ is eventually splitting in $V[\la g_\alpha\ra_{\alpha<\gamma}]$.  By a result of Baumgartner and Dordal~\cite{baumgartner}, $\la x_{j(\alpha)}\ra_{\alpha<\eta}$ is also eventually splitting in $V[G]$.  Choose $\beta<\eta$ such that $x_{j(\alpha)}$ splits $y_G$ for all $\alpha\in\eta\setminus\beta$.  Then there exist $s\in G$ and $\alpha\in\gamma\setminus j(\beta)$ such that $r\geq s\forces \check{\alpha}\in h``J$.  Hence, $\alpha\in I_G$ and $x_\alpha$ splits $y_G$.  Thus, $\la\{x_\alpha,\omega\setminus x_\alpha\}\ra_{\alpha<\kappa}$ witnesses $\cardss_2=\omega_1$ in $V[G]$.
\end{proof}

\begin{definition}
Let $\additivity{\meager}$ denote the additivity of the ideal of meager sets of reals.
\end{definition}

It is known that $\additivity{\meager}\leq\cardb$ and that it is consistent that $\additivity{\meager}<\cardb$.  (See 5.4 and 11.7 of \cite{blass} and 7.3.D of \cite{bartoszynski}).

\begin{corollary}
If $\kappa=\cf\kappa>\omega$, then it is consistent that $\cardss_2=\omega_1$ and $\additivity{\meager}=\cardc=\kappa$.
\end{corollary}
\begin{proof}
Starting with GCH in the ground model, perform a $\kappa$\nbd-long finite support iteration of Hechler forcing.  This forces $\additivity{\meager}=\cardc=\kappa$ (see 11.6 of \cite{blass}).  By Theorem~\ref{THMowhechler}, this also forces $\cardss_2=\omega_1$.
\end{proof}

\section{Models of $\omega_1<\ow{\omega^*}<\cardc$}

To prove the consistency of $\omega_1<\ow{\omega^*}<\cardc$, we employ generalized iteration of forcing along posets as defined by Groszek and Jech~\cite{groszekjech}.  We will only use finite support iterations along well\nbd-founded posets.  For simplicity, we limit our definition of generalized iterations to this special case.

\begin{definition}
Suppose $X$ is a well\nbd-founded poset and $\mbP$ a forcing order consisting of functions on $X$.  Given any $x\in X$, partial map $f$ on $X$, and down\nbd-set $Y$ of $X$, set $\fod{\mbP}{Y}=\{p\restrict Y:p\in\mbP\}$, $\pol{X}{x}=\{y\in X:y<x\}$, $\poleq{X}{x}=\{y\in X:y\leq x\}$, $\fol{\mbP}{x}=\fod{\mbP}{(\pol{X}{x})}$, $\foleq{\mbP}{x}=\fod{\mbP}{(\poleq{X}{x})}$, $\fnl{f}{x}=f\restrict(\pol{X}{x})$, and $\fnleq{f}{x}=f\restrict(\poleq{X}{x})$.  Then $\mbP$ is a \emph{finite support iteration along} $X$ if there exists a sequence $\la\mbQ_x\ra_{x\in X}$ satisfying the following conditions for all $x\in X$ and all $p,q\in\mbP$. 
\begin{enumerate}
\item $\fol{\mbP}{x}$ is a finite support iteration along $\pol{X}{x}$.
\item $\mbQ_x$ is a $(\fol{\mbP}{x})$\nbd-name for a forcing order.
\item $\foleq{\mbP}{x}=\{p\cup\{\la x,q\ra\}:\la p,q\ra\in(\fol{\mbP}{x})*\mbQ_x\}$.
\item $\fnl{\mathbbm{1}_\mbP}{x}\forces \mathbbm{1}_\mbP(x)=\mathbbm{1}_{\mbQ_x}$.
\item $\mbP$ is the set of functions $r$ on $X$ for which $\fnleq{r}{y}\in\foleq{\mbP}{y}$ for all $y\in X$ and $\mathbbm{1}_{\fol{\mbP}{z}}\forces r(z)=\mathbbm{1}_{\mbQ_z}$ for all but finitely many $z\in X$.
\item $p\leq q$ if and only if $\fnl{p}{y}\leq\fnl{q}{y}$ and $\fnl{p}{y}\forces p(y)\leq q(y)$ for all $y\in X$.
\end{enumerate}

Given a finite support iteration $\mbP$ along $X$ and $x\in X$ and a filter $G$ of $\mbP$, set $\gfp{G}{x}=\{p(x):p\in G\}$, $\gfl{G}{x}=\{\fnl{p}{x}:p\in G\}$, and $\gfleq{G}{x}=\{\fnleq{p}{x}:p\in G\}$.  Given any down\nbd-set $Y$ of $X$, set $\gfd{G}{Y}=\{p\restrict Y:p\in G\}$.
\end{definition}

\begin{remark}
If $\mbP$ is a finite support iteration along a well\nbd-founded poset $X$ with down\nbd-set $Y$, then $\fod{\mbP}{Y}$ is an iteration along $Y$, and $\mathbbm{1}_{\fod{\mbP}{Y}}=\mathbbm{1}_\mbP\restrict Y$.
\end{remark}

\begin{definition}
Suppose $\mbP$ is a finite support iteration along a well\nbd-founded poset $X$ with down\nbd-sets $Y$ and $Z$ such that $Y\subseteq Z$.  Then there is a complete embedding $j_Y^Z\colon\fod{\mbP}{Y}\to\fod{\mbP}{Z}$ given by $j_Y^Z(p)=p\cup(\mathbbm{1}_\mbP\restrict Z\setminus Y)$ for all $p\in\fod{\mbP}{Y}$.  This embedding naturally induces an embedding of the class of $(\fod{\mbP}{Y})$\nbd-names, which in turn naturally induces an embedding of the class of atomic forumlae in the $(\fod{\mbP}{Y})$\nbd-forcing language. Let $j_Y^Z$ also denote these embeddings.
\end{definition}

\begin{proposition}\label{PROiterembed}
Suppose $\mbP$, $Y$, and $Z$ are as in the above definition, and $\varphi$ is an atomic formula in the $(\fod{\mbP}{Y})$\nbd-forcing language.  Then, for all $p\in\fod{\mbP}{Z}$, we have $p\forces j_Y^Z(\varphi)$ if and only if $p\restrict Y\forces\varphi$.
\end{proposition}
\begin{proof}
If $p\restrict Y\forces\varphi$, then $p\leq j_Y^Z(p\restrict Y)\forces j_Y^Z(\varphi)$.  Conversely, suppose  $p\restrict Y\not\forces\varphi$.  Then we may choose $q\leq p\restrict Y$ such that $q\forces\neg\varphi$.  Hence, $j_Y^Z(q)\forces\neg j_Y^Z(\varphi)$.  Set $r=q\cup(p\restrict Z\setminus Y)$.  Then $j_Y^Z(q)\geq r\leq p$; hence, $p\not\forces j_Y^Z(\varphi)$.
\end{proof}

\begin{lemma}\label{LEMiterprodlike}
Suppose $\mbP$ is a finite support iteration along a well\nbd-founded poset $X$ and $x$ is a maximal element of $X$.  Set $Y=X\setminus\{x\}$.  Then there is a dense embedding $\phi\colon\mbP\to(\fod{\mbP}{Y})*j_{\pol{X}{x}}^Y(\mbQ_x)$ given by $\phi(p)=\la p\restrict Y,\,j_{\pol{X}{x}}^Y(p(x))\ra$.  Hence, if $G$ is a $\mbP$\nbd-generic filter, then $\gfp{G}{x}$ is $(\mbQ_x)_{\gfl{G}{x}}$\nbd-generic over $V[\gfd{G}{Y}]$.
\end{lemma}
\begin{proof}
First, let us show that $\phi$ is an order embedding. Suppose $r,s\in\mbP$.  Then $r\leq s$ if and only if $r\restrict Y\leq s\restrict Y$ and $\fnl{r}{x}\forces r(x)\leq s(x)$.  Also, $\phi(r)\leq\phi(s)$ if and only if $r\restrict Y\leq s\restrict Y$ and $r\restrict Y\forces j_{\pol{X}{x}}^Y(r(x)\leq s(x))$.  By Proposition~\ref{PROiterembed}, $r\restrict Y\forces j_{\pol{X}{x}}^Y(r(x)\leq s(x))$ if and only if $\fnl{r}{x}\forces r(x)\leq s(x)$; hence, $r\leq s$ if and only if $\phi(r)\leq\phi(s)$.

Finally, let us show that $\ran\phi$ is dense.  Suppose $\la p,q\ra\in(\fod{\mbP}{Y})*j_{\pol{X}{x}}^Y(\mbQ_x)$.  Then there exist $r\leq p$ and $s\in\dom\bigl(j_{\pol{X}{x}}^Y(\mbQ_x)\bigr)$ such that $r\forces s=q\in j_{\pol{X}{x}}^Y(\mbQ_x)$.  Hence, $\la r,s\ra\leq\la p,q\ra$.  Also, $s$ is a $(j_{\pol{X}{x}}^Y``(\fol{\mbP}{x}))$\nbd-name; hence, there exists a $(\fol{\mbP}{x})$\nbd-name $t$ such that $j_{\pol{X}{x}}^Y(t)=s$.  Hence, $r\forces j_{\pol{X}{x}}^Y(t\in\mbQ_x)$; hence, $\fnl{r}{x}\forces t\in\mbQ_x$.  Hence, $r\cup\{\la x,t\ra\}\in\mbP$ and $\phi(r\cup\{\la x,t\ra\})=\la r,s\ra$.  Thus, $\ran\phi$ is dense.
\end{proof}

\begin{remark}
Proposition~\ref{PROiterembed} and Lemma~\ref{LEMiterprodlike} and their proofs remain valid for arbitrary iterations along posets as defined in \cite{groszekjech}.
\end{remark}

\begin{lemma}\label{LEMgenericpseudointersection}
Let $\mbP$ be a forcing order, $A$ a subset of $[\omega]^\omega$ with the SFIP, $\mbQ$ the Booth forcing for $A$, $x$ a $\mbQ$\nbd-name for a generic pseudointersection of $A$, and $B$ a $\mbP$\nbd-name such that $\mathbbm{1}_{\mbP}$ forces $\check{A}\subseteq B\subseteq[\omega]^\omega$ and forces $B$ to have the SFIP.  Let $i$ and $j$ be the canonical embeddings, respectivly, of $\mbP$\nbd-names and $\mbQ$\nbd-names into $(\mbP*\check{\mbQ})$\nbd-names.  Then $\mathbbm{1}_{\mbP*\check{\mbQ}}$ forces $i(B)\cup\{j(x)\}$ to have the SFIP.
\end{lemma}
\begin{proof}
Seeking a contradiction, suppose $r_0=\la p_0,\la\sigma,F\ra^{\check{}}\ra\in\mbP*\check{\mbQ}$ and $n<\omega$ and $p_0\forces H\in[B]^{<\omega}$ and $r_0\forces j(x)\cap\bigcap i(H)\subseteq\check{n}$.  Then $p_0$ forces $\check{F}\cup H\subseteq B$, which is forced to have the SFIP; hence, there exist $p_1\leq p_0$ and $m\in\omega\setminus n$ such that $p_1\forces\check{m}\in\bigcap(\check{F}\cup H)$.  Set $r_1=\la p_1,\la\sigma\cup\{m\},F\ra^{\check{}}\ra$.  Then $r_0\geq r_1\forces \check{m}\in j(x)\cap\bigcap i(H)$, contradicting how we chose $r_0$.
\end{proof}

\begin{lemma}\label{LEMpropKcccabs}
Suppose $\mbP$ and $\mbQ$ are forcing orders such that $\mbP$ is ccc and $\mbQ$ has property (K).  Then $\mathbbm{1}_\mbP$ forces $\check{\mbQ}$ to have property (K).
\end{lemma}
\begin{proof}
Suppose the lemma fails.  Then there exist $p\in\mbP$ and $f$ such that $p\forces f\in\check{\mbQ}^{\omega_1}$ and $p\forces\forall J\in[\omega_1]^{\omega_1}\ \exists \alpha,\beta\in J\ f(\alpha)\perp f(\beta)$.  For each $\alpha<\omega_1$, choose $p_\alpha\leq p$ and $q_\alpha\in\mbQ$ such that $p_\alpha\forces f(\alpha)=\check{q}_\alpha$.  Then there exists $I\in[\omega_1]^{\omega_1}$ such that $q_\alpha\not\perp q_\beta$ for all $\alpha,\beta\in I$.  Let $J$ be the $\mbP$\nbd-name $\{\la\check{\alpha},p_\alpha\ra:\alpha\in I\}$.  Then $p\forces\forall\alpha,\beta\in J\ f(\alpha)=\check{q}_\alpha\not\perp\check{q}_\beta=f(\beta)$.  Hence, $p\forces\card{J}\leq\omega$.  Since $\mbP$ is ccc, there exists $\alpha\in I$ such that $p\forces J\subseteq\check{\alpha}$.  But this contradicts $p\geq p_\alpha\forces\check{\alpha}\in J$.
\end{proof}

\begin{lemma}\label{LEMpropK}
Suppose $\mbP$ is a finite support iteration along a well\nbd-founded poset $X$ and $\fnl{\mathbbm{1}_\mbP}{x}$ forces $\mbQ_x$ to have property (K) for all $x\in X$.  Then $\mbP$ has property (K).
\end{lemma}
\begin{proof}
\Wma\ the lemma holds whenever $X$ is replaced by a poset of lesser height.  Let $I\in[\mbP]^{\omega_1}$.  \Wma\ $\{\supp{p}:p\in I\}$ is a $\Delta$\nbd-system; let $\sigma$ be its root.  Set $Y_0=\bigcup_{x\in\sigma}\pol{X}{x}$.  Then $\fod{\mbP}{Y_0}$ has property (K).  Let $n=\card{\sigma\setminus Y_0}$ and $\la x_i\ra_{i<n}$ biject from $n$ to $\sigma\setminus Y_0$.  Set $Y_{i+1}=Y_i\cup\{x_i\}$ for all $i<n$.  Suppose $i<n$ and $\fod{\mbP}{Y_i}$ has property (K).  By Lemma~\ref{LEMpropKcccabs}, $\mathbbm{1}_{\fod{\mbP}{Y_i}}$ forces $j_{\pol{X}{x_i}}^{Y_i}(\mbQ_{x_i})$ to have property (K).  Hence, $\fod{\mbP}{Y_{i+1}}$ has property (K), for it densely embeds into $\fod{\mbP}{Y_i}*j_{\pol{X}{x_i}}^{Y_i}(\mbQ_{x_i})$ by Lemma~\ref{LEMiterprodlike}.  By induction, $\fod{\mbP}{Y_n}$ has property (K); hence, there exists $J\in[I]^{\omega_1}$ such that $p\restrict Y_n\not\perp q\restrict Y_n$ for all $p,q\in J$.  Fix $p,q\in J$ and choose $r$ such that $r\leq p\restrict Y_n$ and $r\leq q\restrict Y_n$.  Set $s=r\cup(p\restrict\supp{p}\setminus Y_n)\cup (q\restrict\supp{q}\setminus Y_n)$ and $t=s\cup(\mathbbm{1}_{\mbP}\restrict X\setminus\dom s)$.  Then $t\leq p,q$.
\end{proof}

\begin{lemma}\label{LEMkappamuposet}
Suppose $\cf\kappa=\kappa\leq\lambda=\lambda^{<\kappa}$.  Then there exists a $\kappa$\nbd-like, $\kappa$\nbd-directed, well\nbd-founded poset $\Xi$ with cofinality and cardinality $\lambda$.
\end{lemma}
\begin{proof}
Let $\{x_\alpha:\alpha<\lambda\}$ biject from $\lambda$ to $[\lambda]^{<\kappa}$.  Construct $\la y_\alpha\ra_{\alpha<\lambda}\in([\lambda]^{<\kappa})^\lambda$ as follows.  Given $\alpha<\lambda$ and $\la y_\beta\ra_{\beta<\alpha}$, choose $\xi_\alpha\in\lambda\setminus\bigcup_{\beta<\alpha}y_\beta$ and set $y_\alpha=x_\alpha\cup\{\xi_\alpha\}$.  Let $\Xi$ be $\{y_\alpha:\alpha<\lambda\}$ ordered by inclusion.  Then $\Xi$ is cofinal with $[\lambda]^{<\kappa}$; hence, $\Xi$ is $\kappa$\nbd-directed and has cofinality $\lambda$.  Also, $\Xi$ is well\nbd-founded because $\la y_\alpha\ra_{\alpha<\lambda}$ is nondecreasing.  Finally, $\Xi$ is $\kappa$\nbd-like because for all $I\in[\lambda]^\kappa$ we have $\card{\bigcup_{\alpha\in I}y_\alpha}\geq\card{\{\xi_\alpha:\alpha\in I\}}=\kappa$; whence, $\{y_\alpha:\alpha\in I\}$ has no upper bound in $[\lambda]^{<\kappa}$.
\end{proof}

\begin{definition}
A point $q$ in a space $X$ is a \emph{$P_\kappa$\nbd-point} if every intersection of fewer than $\kappa$\nbd-many neighborhoods of $q$ contains a neighborhood of $q$.
\end{definition}

\begin{definition}
For all $x,y\subseteq\omega$, define $x\subseteq^*y$ as $\card{x\setminus y}<\omega$.  Let $\cardp$ denote the minimum value of $\card{A}$ where $A$ ranges over the subsets of $[\omega]^\omega$ that have SFIP yet have no pseudointersection.
\end{definition}

\begin{remark}
It easily seen that $\omega_1\leq\cardp\leq\cardh$.
\end{remark}

\begin{theorem}\label{THMowomegastarkappa}
Suppose $\omega_1\leq\cf\kappa=\kappa\leq\lambda=\lambda^{<\kappa}$.  Then there is a property (K) forcing extension in which
\begin{equation*}
\cardp=\opi{\omega^*}=\ow{\omega^*}=\cardss_2=\cardb=\kappa\leq\lambda=\cardc.
\end{equation*}
Moreover, in this extension $\omega^*$ has $P_\kappa$\nbd-points; whence, $\max_{q\in\omega^*}\ochar{q,\omega^*}=\kappa$.
\end{theorem}
\begin{proof}
Let $\Xi$ be as in Lemma~\ref{LEMkappamuposet}.  Let $\la \sigma_\alpha\ra_{\alpha<\lambda}$ biject from $\lambda$ to $\Xi$.  Let $\la\la\zeta_\alpha,\eta_\alpha\ra\ra_{\alpha<\lambda}$ biject from $\lambda$ to $\lambda^2$.  Given $\alpha<\lambda$ and $\la\tau_{\zeta_\beta,\eta_\beta}\ra_{\beta<\alpha}\in\Xi^\alpha$, choose $\tau_{\zeta_\alpha,\eta_\alpha}\in\Xi$ such that $\sigma_{\zeta_\alpha}<\tau_{\zeta_\alpha,\eta_\alpha}\not\leq\tau_{\zeta_\beta,\eta_\beta}$ for all $\beta<\alpha$.  We may so choose $\tau_{\zeta_\alpha,\eta_\alpha}$ because $\Xi$ is directed and has cofinality $\lambda$.

Let us construct a finite support iteration $\mbP$ along $\Xi$.  Since $\Xi$ is well\nbd-founded, we may define $\mbQ_\sigma$ in terms of $\fol{\mbP}{\sigma}$ for each $\sigma\in\Xi$.  Suppose $\sigma\in\Xi$ and, for all $\tau<\sigma$, we have $\card{\foleq{\mbP}{\tau}}<\kappa$ and $\mathbbm{1}_{\fol{\mbP}{\tau}}$ forces $\mbQ_\tau$ to have property (K).  Then $\fol{\mbP}{\sigma}$ has property (K) by Lemma~\ref{LEMpropK}, and hence is ccc.  Moreover, $\card{\fol{\mbP}{\sigma}}<\kappa$ because $\fol{\mbP}{\sigma}$ is a finite support iteration along $\pol{\Xi}{\sigma}$ and $\card{\pol{\Xi}{\sigma}}<\kappa$.  Hence, $\mathbbm{1}_{\fol{\mbP}{\sigma}}\forces \card{\cardc^{<\kappa}}\leq((\kappa^\omega)^{<\kappa})^{\check{}}\leq\lambda$.  Let $\mcE_\sigma$ be a $(\fol{\mbP}{\sigma})$\nbd-name for the set of all $E$ in the $(\fol{\mbP}{\sigma})$\nbd-generic extension for which $E\in[[\omega]^\omega]^{<\kappa}$ and $E$ has the SFIP.  Then we may choose a $(\fol{\mbP}{\sigma})$\nbd-name $f_\sigma$ such that $\mathbbm{1}_{\fol{\mbP}{\sigma}}$ forces $f_\sigma$ to be a surjection from $\lambda$ to $\mcE_\sigma$.  \Wma\ we are given corresponding $f_\tau$ for all $\tau<\sigma$.  If there exist $\alpha,\beta<\lambda$ such that $\sigma=\tau_{\alpha,\beta}$, then let $\mbQ_\sigma$ be a $(\fol{\mbP}{\sigma})$\nbd-name for $\mbQ_\sigma'\times\Fn{\omega}{2}$ where $\mbQ_\sigma'$ is a $(\fol{\mbP}{\sigma})$\nbd-name for the Booth forcing for $f_{\sigma_\alpha}(\beta)$.  If there are no such $\alpha$ and $\beta$, then let $\mbQ_\sigma$ be a $(\fol{\mbP}{\sigma})$\nbd-name for a singleton poset.  Then $\mathbbm{1}_{\fol{\mbP}{\sigma}}$ forces $\mbQ_\sigma$ to have property (K).  Also, \wma\ $\card{\mbQ_\sigma}<\kappa$.  Hence, $\card{\foleq{\mbP}{\sigma}}<\kappa$.

By induction, $\card{\foleq{\mbP}{\sigma}}<\kappa$ and $\mathbbm{1}_{\fol{\mbP}{\sigma}}$ forces $\mbQ_\sigma$ to have property (K) for all $\sigma\in\Xi$.  Hence, $\mbP$ has property (K) by Lemma~\ref{LEMpropK}, and hence is ccc.  Also, since $\card{\Xi}\leq\lambda$ and $\mbP$ is a finite support iteration, $\card{\mbP}\leq\lambda$.  Let $G$ be a $\mbP$\nbd-generic filter.  Then $\cardc^{V[G]}\leq\lambda^\omega=\lambda$.  Moreover, $\cardc^{V[G]}\geq\lambda$ because $\mbP$ adds $\lambda$\nbd-many Cohen reals.

By Theorem~\ref{THMopicardh} and Lemma~\ref{LEMowsplitterupperbound}, it suffices to show that $\cardb^{V[G]}\leq\kappa\leq\cardp^{V[G]}$, that $\cardss_2^{V[G]}\leq\kappa$, and that some $q\in(\omega^*)^{V[G]}$ is a $P_\kappa$\nbd-point.  First, we prove $\kappa\leq\cardp^{V[G]}$.  Suppose $E\in([[\omega]^\omega]^{<\kappa})^{V[G]}$ and $E$ has the SFIP.  Then there exists $\alpha<\lambda$ such that $E\in V[\gfl{G}{\sigma_\alpha}]$ because $\Xi$ is $\kappa$\nbd-directed.  Hence, there exists $\beta<\lambda$ such that $(f_{\sigma_\alpha})_{\gfl{G}{\sigma_\alpha}}(\beta)=E$.  Hence, $E$ has a pseudointersection in $V[\gfleq{G}{\tau_{\alpha,\beta}}]$.  Thus, $\kappa\leq\cardp^{V[G]}$.

Second, let us show that $\cardb^{V[G]}\leq\kappa$.  For each $\alpha<\kappa$, let $u_\alpha$ be the increasing enumeration of the Cohen real added by the $\Fn{\omega}{2}$ factor of $\mbQ_{\tau_{0,\alpha}}$.  Then it suffices to show that $\{u_\alpha:\alpha<\kappa\}$ is unbounded in $(\omega^\omega)^{V[G]}$.  Suppose $v\in(\omega^\omega)^{V[G]}$.  Then there exists $\sigma\in\Xi$ such that $v\in V[\gfl{G}{\sigma}]$.  Since $\Xi$ is $\kappa$\nbd-like, there exists $\alpha<\kappa$ such that $\tau_{0,\alpha}\not\leq\sigma$.  By Lemma~\ref{LEMiterprodlike}, $u_\alpha$ enumerates a real Cohen generic over $V[\gfl{G}{\sigma}]$; hence, $u_\alpha$ is not eventually dominated by $v$.

Third, let us prove $\cardss_2^{V[G]}\leq\kappa$.  For each $\alpha<\lambda$, let $x_\alpha$ be the Cohen real added by the $\Fn{\omega}{2}$ factor of $\mbQ_{\tau_{0,\alpha}}$.  Suppose $I\in([\lambda]^\kappa)^{V[G]}$ and $y\in([\omega]^\omega)^{V[G]}$.   Then there exists $\sigma\in\Xi$ such that $y\in V[\gfl{G}{\sigma}]$.  Since $\Xi$ is $\kappa$\nbd-like, there exists $\alpha\in I$ such that $\tau_{0,\alpha}\not\leq\sigma$.  By Lemma~\ref{LEMiterprodlike}, $x_\alpha$ is Cohen generic over $V[\gfl{G}{\sigma}]$, and therefore splits $y$.  Thus, $\la\{x_\alpha,\omega\setminus x_\alpha\}\ra_{\alpha<\lambda}$ witnesses $\cardss_2^{V[G]}\leq\kappa$.

Finally, let us construct a $P_\kappa$\nbd-point $q\in(\omega^*)^{V[G]}$.  Let $\sqsubseteq$ be an extension of the ordering of $\Xi$ to a well\nbd-ordering of $\Xi$.  For each $\sigma\in\Xi$, set $Y_\sigma=\{\tau\in\Xi:\tau\sqsubset\sigma\}$.  Set $\rho=\min_\sqsubseteq\Xi$ and choose $U_\rho\in(\omega^*)^V$.  Suppose $\tau\in\Xi$ and $\sigma$ is a final predecessor of $\tau$ with respect to $\sqsubseteq$ and $U_\sigma\in(\omega^*)^{V[\gfd{G}{Y_\sigma]}}$.  If there are no $\alpha,\beta<\lambda$ such that $\sigma=\tau_{\alpha,\beta}$ and $(f_{\sigma_\alpha})_{\gfl{G}{\sigma_\alpha}}(\beta)\subseteq U_\sigma$, then choose $U_\tau\in(\omega^*)^{V[\gfd{G}{Y_\tau]}}$ such that $U_\tau\supseteq U_\sigma$.  Now suppose such $\alpha$ and $\beta$ exist.  Let $v_\sigma$ be the pseudointersection of $(f_{\sigma_\alpha})_{\gfl{G}{\sigma_\alpha}}(\beta)$ added by $\mbQ_{\sigma}'$.  

By Lemmas~\ref{LEMiterprodlike} and \ref{LEMgenericpseudointersection}, $U_\sigma\cup\{v_\sigma\}$ has the SFIP; hence, we may choose $U_\tau\in(\omega^*)^{V[\gfd{G}{Y_\tau]}}$ such that $U_\tau\supseteq U_\sigma\cup\{v_\sigma\}$.  For $\tau\in\Xi$ that are limit points with respect to $\sqsubseteq$, choose $U_\tau\in(\omega^*)^{V[\gfd{G}{Y_\tau}]}$ such that $U_\tau\supseteq\bigcup_{\sigma\sqsubset\tau}U_\sigma$; set $q=\bigcup_{\tau\in\Xi}U_\tau$.  Then, arguing as in the proof of $\kappa\leq\cardp^{V[G]}$, we have that $q$ is a $P_\kappa$\nbd-point in $(\omega^*)^{V[G]}$.
\end{proof}
  
The forcing extension of Theorem~\ref{THMowomegastarkappa} can be modified to satisfy $\cardb=\cards<\ow{\omega^*}<\cardc$.

\begin{definition}
Given a class $\mcJ$ of posets and a cardinal $\kappa$, let $\MA{\kappa}{\mcJ}$ denote the statement that, given any $\mbP\in\mcJ$ and fewer than $\kappa$\nbd-many dense subsets of $\mbP$, there is a filter of $\mbP$ intersecting each of these dense sets.  We may replace $\mcJ$ with a descriptive term for $\mcJ$ when there is no ambiguity.  For example, $\MA{\cardc}{\text{ccc}}$ is Martin's axiom.
\end{definition}

\begin{theorem}\label{THMcardblessowomegastar}
Suppose $\omega_1<\cf\kappa=\kappa\leq\lambda=\lambda^{<\kappa}$.  Then there is a property (K) forcing extension in which
\begin{equation*}
\omega_1=\opi{\omega^*}=\cardb=\cards<\ow{\omega^*}=\cardss_2=\kappa\leq\lambda=\cardc.
\end{equation*}
\end{theorem}
\begin{proof}
Let $\mbP$ be as in the proof of Theorem~\ref{THMowomegastarkappa}.  Set $\mbR=\mbP\times\Fn{\omega_1}{2}$, which has property (K) because $\mbP$ does.  Let $K$ be a generic filter of $\mbR$. Let $\pi_0$ and $\pi_1$ be the natural coordinate projections on $\mbR$; let $\pi_0$ and $\pi_1$ also denote their respective natural extensions to the class of $\mbR$\nbd-names.  Set $G=\pi_0``K$ and $H=\pi_1``K$.  Then $\cardc^{V[K]}=\lambda$ clearly holds.  Adding $\omega_1$\nbd-many Cohen reals to any model of ZFC forces $\cardb=\cards=\omega_1$, and $\opi{\omega^*}=\cardh\leq\cardb$, so $\opi{\omega^*}^{V[K]}=\cardb^{V[K]}=\cards^{V[K]}=\omega_1$.  

For each $\alpha<\lambda$, let $x_\alpha$ be the Cohen real added by the $\Fn{\omega}{2}$ factor of $\mbQ_{\tau_{0,\alpha}}$.  Suppose $I\in([\lambda]^\kappa)^{V[K]}$ and $y\in([\omega]^\omega)^{V[K]}$.   Then there exists $\sigma\in\Xi$ such that $y\in V[(\gfl{G}{\sigma})\times H]$.  Since $\Xi$ is $\kappa$\nbd-like, there exists $\alpha\in I$ such that $\tau_{0,\alpha}\not\leq\sigma$.  By Lemma~\ref{LEMiterprodlike}, $x_\alpha$ is Cohen generic over $V[\gfl{G}{\sigma}]$; hence, $x_\alpha$ is Cohen generic over $V[(\gfl{G}{\sigma})\times H]$ and therefore splits $y$.  Thus, $\la\{x_\alpha,\omega\setminus x_\alpha\}\ra_{\alpha<\lambda}$witnesses $\cardss_2^{V[K]}\leq\kappa$.

Therefore, it suffices to show that $\ow{\omega^*}^{V[K]}\geq\kappa$.  Suppose $\mu<\kappa$ and $\mcA$ is an $\mbR$\nbd-name for a base of $\omega^*$.  Choose an $\mbR$\nbd-name $q$ for an element of $\omega^*$ with character $\lambda$.  Let $f$ be a name for an injection from $\lambda$ into $\mcA$ such that $q\in\bigcap\ran f$.  Let $g$ be a name for an element of $([\omega]^\omega)^\lambda$ such that $q\in g(\alpha)^*\subseteq f(\alpha)$ for all $\alpha<\lambda$.  For each $\alpha<\lambda$, let $u_\alpha$ be a name for $g(\alpha)$ such that $u_\alpha=\{\{\check{n}\}\times A_{\alpha,n}:n<\omega\}$ where each $A_{\alpha,n}$ is a countable antichain of $\mbR$.  Since $\max\{\omega_1,\mu\}<\lambda$, there exist $\xi<\omega_1$ and $J\in[\lambda]^\mu$ such that $\ran\pi_1(u_\alpha)\subseteq \Fn{\xi}{2}$ for all $\alpha\in J$.  It suffices to show that $\{(u_\alpha)_K:\alpha\in J\}$ has a pseudointersection in $V[K]$.

For each $\alpha\in J$, set $v_\alpha=\{\la \check{n},r\ra:\la \check{n},\la p,r\ra\ra\in u_\alpha\text{ and }p\in G\}$.  Set $H_0=H\cap\Fn{\xi}{2}$.  By Bell's Theorem~\cite{bell}, $\MA{\cardp}{\sigma\text{-centered}}$ is a theorem of ZFC.  Hence, $V[G]$ satisfies $\MA{\kappa}{\sigma\text{-centered}}$.  By an argument of Baumgartner and Tall communicated by Roitman~\cite{roitman}, adding a single Cohen real preserves $\MA{\kappa}{\sigma\text{-centered}}$.  Since Booth forcing for $\{(v_\alpha)_{H_0}:\alpha\in J\}$ is $\sigma$\nbd-centered, $\{(v_\alpha)_{H_0}:\alpha\in J\}$, which is equal to $\{(u_\alpha)_K:\alpha\in J\}$, has a pseudointersection in $V[G\times H_0]$.
\end{proof}

\section{Local Noetherian type and $\pi$\nbd-type}

\begin{definition}
For every infinite cardinal $\kappa$, let $u(\kappa)$ denote the space of uniform ultrafilters on $\kappa$.
\end{definition}

Dow and Zhou~\cite{dow} proved that there is a point in $\omega^*$ that (along with satisfying some additional properties) has an \omegaop\ local base.  We present a simpler construction of an \omegaop\ local base which also naturally generalizes to every $u(\kappa)$.  This construction is essentially due to Isbell~\cite{isbell}, who was interested in actual intersections as opposed to pseudointersections.

\begin{definition}
Given cardinals $\lambda\geq\kappa\geq\omega$ and a point $p$ in a space $X$, a \emph{local} $\la\lambda,\kappa\ra$\nbd-\emph{splitter} is a set $\mcU$ of $\lambda$\nbd-many open neighborhoods of $p$ such that $p$ is not in the interior of $\bigcap\mcV$ for any $\mcV\in[\mcU]^\kappa$.
\end{definition}

\begin{lemma}\label{LEMkappasizeposet}
Every poset $P$ is almost \cardop{\card{P}}.
\end{lemma}
\begin{proof}
Let $\kappa=\card{P}$ and let $\la p_\alpha\ra_{\alpha<\kappa}$ biject from $\kappa$ to $P$.  Define a partial map $f\colon\kappa\to P$ as follows.  Suppose $\alpha<\kappa$ and we have a partial map $f_\alpha\colon\alpha\to P$.  If $\ran f_\alpha$ is dense in $P$, then set $f_{\alpha+1}=f_\alpha$.  Otherwise, set $\beta=\min\{\delta<\kappa:p_\delta\not\geq q\text{ for all }q\in\ran f_\alpha\}$ and set $f_{\alpha+1}=f_\alpha\cup\{\la\alpha,p_\beta\ra\}$.  For limit ordinals $\gamma\leq\kappa$, set $f_\gamma=\bigcup_{\alpha<\gamma}f_\alpha$.  Set $f=f_\kappa$.  Then $f$ is nonincreasing; hence, $\ran f$ is \kappaop.  Moreover, $\ran f$ is dense in $P$.
\end{proof}

\begin{lemma}\label{LEMocharconverse}
Suppose $X$ is a space with a point $p$ at which there is no finite local base.  Then $\ochar{p,X}$ is the least $\kappa\geq\omega$ for which there is a local $\la\character{p,X},\kappa\ra$\nbd-splitter at $p$.  Moreover, if $\lambda>\character{p,X}$, then $p$ does not have a local $\la\lambda,\kappa\ra$\nbd-splitter at $p$ for any $\kappa<\lambda$ or $\kappa\leq\cf\lambda$.
\end{lemma}
\begin{proof}
By Lemma~\ref{LEMkappasizeposet}, $\character{p,X}\geq\ochar{p,X}$; hence, a \cardop{\ochar{p,X}}\ local base at $p$ (which necessarily has size $\character{p,X}$) is a local $\la\character{p,X},\ochar{p,X}\ra$\nbd-splitter at $p$.  To show the converse, let $\lambda=\character{p,X}$ and let $\la U_\alpha\ra_{\alpha<\lambda}$ be a sequence of open neighborhoods of $p$.  Let $\{V_\alpha:\alpha<\lambda\}$ be a local base at $p$.  For each $\alpha<\lambda$, choose $W_\alpha\in\{V_\beta:\beta<\lambda\}$ such that $W_\alpha\subseteq U_\alpha\cap V_\alpha$.  Then $\{W_\alpha:\alpha<\lambda\}$ is a local base at $p$.  Let $\kappa<\ochar{p,X}$.  Then there exist $\alpha<\lambda$ and $I\in[\lambda]^\kappa$ such that $W_\alpha\subseteq\bigcap_{\beta\in I}W_\beta$.  Hence, $p$ is in the interior of $\bigcap_{\beta\in I}U_\beta$.  Hence, $\{U_\alpha:\alpha<\lambda\}$ is not a local $\la\lambda,\kappa\ra$\nbd-splitter at $p$.

To prove the second half of the lemma, suppose $\lambda>\character{p,X}$ and $\mcA$ is a set of $\lambda$\nbd-many open neighborhoods of $p$.  Let $\mcB$ be a local base at $p$ of size $\character{p,X}$.  Then, for all $\kappa<\lambda$ and $\kappa\leq\cf\lambda$, there exist $U\in\mcB$ and $\mcC\in[\mcA]^\kappa$ such that $U\subseteq\bigcap\mcC$.  Hence, $\mcA$ is not a local $\la\lambda,\kappa\ra$\nbd-splitter at $p$.
\end{proof}

\begin{theorem}\label{THMocharomegainukappa}
For each $\kappa\geq\omega$, there exists $p\in u(\kappa)$ such that $\ochar{p,u(\kappa)}=\omega$ and $\character{p,u(\kappa)}=2^\kappa$.
\end{theorem}
\begin{proof}
Let $A$ be an independent family of subsets of $\kappa$ of size $2^\kappa$.  Set $B=\bigcup_{F\in[A]^\omega}\{x\subseteq\kappa:\forall y\in F\ \ \card{x\setminus y}<\kappa\}$.  Since $A$ is independent, we may extend $A$ to an ultrafilter $p$ on $\kappa$ such that $p\cap B=\emptyset$.  For each $x\subseteq\kappa$, set $x^*=\{q\in u(\kappa):x\in q\}$.  Then  $\{x^*:x\in A\}$ is a local $\la 2^\kappa,\omega\ra$\nbd-splitter at $p$.  Since $\character{p,u(\kappa)}\leq 2^\kappa$, it follows from Lemma~\ref{LEMocharconverse} that $\ochar{p,u(\kappa)}=\omega$ and $\character{p,u(\kappa)}=2^\kappa$.
\end{proof}

\begin{definition}
Let $\carda$ denote the minimum of the cardinalities of infinite, maximal almost disjoint subfamilies of $[\omega]^\omega$.  Let $\cardi$ denote the minimum of the cardinalities of infinite, maximal independent subfamilies of $[\omega]^\omega$.
\end{definition}

It is known that $\cardb\leq\carda$ and $\cardr\leq\cardi\geq\cardd\geq\cards$.  (See 8.4, 8.12, 8.13 and 3.3 of \cite{blass}.)  Because of Kunen's result that $\carda=\aleph_1$ in the Cohen model (see VIII.2.3 of \cite{kunen}), it is consistent that $\carda<\cardr$.  Also, Shelah~\cite{shelahad} has constructed a model of $\cardr\leq\cardu<\carda$.

In ZFC, the best upper bound of $\ochar{\omega^*}$ of which we know is $\cardc$ by Lemma~\ref{LEMkappasizeposet}.  We will next prove Theorem~\ref{THMlocntabovei}, which implies that, except for $\cardc$ and possibly $\cf\cardc$, all of the cardinal characteristics of the continuum with definitions included in Blass~\cite{blass} can consistently be simultaneously strictly less than $\ochar{\omega^*}$.

\begin{lemma}\label{LEMkappacrosslambda}
Suppose $\kappa$, $\lambda$, and $\mu$ are regular cardinals and $\kappa\leq\lambda>\mu$.  Then $(\kappa\times\lambda)^{\mathrm{op}}$ is not almost \cardop{\mu}.
\end{lemma}
\begin{proof}
Let $I$ be a cofinal subset of $\kappa\times\lambda$.  Then it suffices to show that $I$ is not $\mu$\nbd-like.  If $\kappa=\lambda$, then $I$ is not $\mu$\nbd-like because it is $\lambda$\nbd-directed.  Suppose $\kappa<\lambda$.  Then there exists $\alpha<\kappa$ such that $\card{I\cap(\{\alpha\}\times\lambda)}=\lambda$; hence, $I$ has an increasing $\lambda$\nbd-sequence; hence, $I$ is not $\mu$\nbd-like.
\end{proof}

\begin{lemma}\label{LEMmif}
Given any infinite independent subfamily $I$ of $[\omega]^\omega$, there exists $J\subseteq[\omega]^\omega$ such that if $x$ is a generic pseudointersection of $J$ then $I\cup\{x\}$ is independent, but $I\cup\{x,y\}$ is not independent for any $y\in[\omega]^\omega\cap V\setminus I$.
\end{lemma}
\begin{proof}
See Exercise A12 on page 289 of Kunen~\cite{kunen}.
\end{proof}

\begin{definition}
We say a $P_\kappa$\nbd-point in a space is \emph{simple} if it has a local base of order type $\kappa^{\mathrm{op}}$.
\end{definition}

\begin{theorem}\label{THMlocntabovei}
Suppose $\omega_1\leq\cf\kappa=\kappa\leq\cf\lambda=\lambda=\lambda^{<\kappa}$.  Then there is a property (K) forcing extension satisfying $\cardp=\carda=\cardi=\cardu=\kappa\leq\lambda=\ochar{\omega^*}=\cardc$.
\end{theorem}
\begin{proof}
We will construct a finite support iteration $\la\mbP_\alpha\ra_{\alpha\leq\lambda\kappa}$ where $\lambda\kappa$ denotes the ordinal product of $\lambda$ and $\kappa$.  It suffices to ensure that the iteration is at every stage property (K) and of size at most $\lambda$, and that $V^{\mbP_{\lambda\kappa}}$ satisfies $\max\{\carda,\cardi,\cardu\}\leq\kappa\leq\cardp$ and $\lambda\leq\ochar{\omega^*}$.  Our strategy is to interleave an iteration of length $\lambda\kappa$ and three iterations of length $\kappa$.  At every stage below $\lambda\kappa$, add another piece of what will be an ultrafilter base that, ordered by $\supseteq^*$, will be isomorphic to a cofinal subset of $\kappa\times\lambda$.  Also, at every stage we will add a pseudointersection, such that the final model satisfies $\cardp\geq\kappa$.  After each limit stage of cofinality $\lambda$, add an element to each of three objects that, when completed, will be a maximal almost disjoint family of size $\kappa$, a maximal independent family of size $\kappa$, and a base of a simple $P_\kappa$\nbd-point in $\omega^*$.

Let $\varphi\colon\lambda^2\rightarrow\lambda$ be a bijection such that $\varphi(\alpha,\beta)\geq\alpha$ for all $\alpha,\beta<\lambda$.  For each $\la\alpha,\beta\ra\in\kappa\times\lambda$, set $E_{\alpha,\beta}=\{\la\gamma,\delta\ra\in\kappa\times\lambda:\lambda\gamma+\delta<\lambda\alpha+\beta\}$.  Suppose $\la\alpha,\beta\ra\in\kappa\times\lambda$ and we have constructed $\la\mbP_\gamma\ra_{\gamma\leq\lambda\alpha+\beta}$ to have property (K) and size at most $\lambda$ at all of its stages, and a sequence $\la x_{\gamma,\delta}\ra_{\la\gamma,\delta\ra\in E_{\alpha,\beta}}$ of $\mbP_{\lambda\alpha+\beta}$\nbd-names each forced to be in $[\omega]^{\omega}$.  Set $B=\{x_{\gamma,\delta}:\la\gamma,\delta\ra\in E_{\alpha,\beta}\}$.  Let $\la S_\gamma\ra_{\gamma<\kappa}$ be a partition of $\lambda$ into $\kappa$\nbd-many stationary sets such that $S_0$ contains all successor ordinals.  Suppose we have constructed a sequence $\la\rho_{\gamma,\delta}\ra_{\la\gamma,\delta\ra\in E_{\alpha,\beta}}\in\lambda^{E_{\alpha,\beta}}$ such that we always have $\rho_{\gamma,\delta}\in S_\gamma$ and $\rho_{\gamma,\delta_0}<\rho_{\gamma,\delta_1}$ whenever $\delta_0<\delta_1$.  Set $D_{\alpha,\beta}=\{\la\gamma,\rho_{\gamma,\delta}\ra:\la\gamma,\delta\ra\in E_{\alpha,\beta}\}$.  Further suppose that $\{\la\la\gamma,\rho_{\gamma,\delta}\ra,x_{\gamma,\delta}\ra:\la\gamma,\delta\ra\in E_{\alpha,\beta}\}$ is forced to be an order embedding of $D_{\alpha,\beta}$ into $\la[\omega]^\omega,\supseteq^*\ra$ and that its range $B$ is forced to have the SFIP.  Also suppose that we have the following if $\alpha>0$.
\begin{equation}\label{eqlattice}
\forces_{\lambda\alpha+\beta}\forall\sigma\in[B]^{<\omega}\ \exists\delta<\lambda\ \,\bigcap\sigma\not\subseteq^*x_{0,\delta}
\end{equation}
For each $\varepsilon<\lambda$, set $A_\varepsilon=\{x_{\gamma,\delta}:\la\gamma,\delta\ra\in E_{\alpha,\beta}\text{ and }\la\gamma,\rho_{\gamma,\delta}\ra<\la\alpha,\varepsilon\ra\}$. 

Let $y_\beta$ be a $\mbP_{\lambda\alpha+\beta}$\nbd-name for a surjection from $\lambda$ to $[\omega]^\omega$.  \Wma\ that corresponding $y_\gamma$ have already been constructed for all $\gamma<\beta$.  Let $\varphi(\zeta,\eta)=\beta$.  

\begin{claim}
If $\alpha>0$, then we may choose $z\in\{y_\zeta(\eta),\,\omega\setminus y_\zeta(\eta)\}$ such that
\begin{equation*}
\forces_{\lambda\alpha+\beta}\forall\sigma\in[B]^{<\omega}\ \exists\delta<\lambda\ \,z\cap\bigcap\sigma\not\subseteq^*x_{0,\delta}.
\end{equation*}
\end{claim}
\begin{proof}
Suppose not.  Let $\{z_0,z_1\}=\{y_\zeta(\eta),\,\omega\setminus y_\zeta(\eta)\}$.  Then, working in a generic extension by $\mbP_{\lambda\alpha+\beta}$, there exist $\sigma_0,\sigma_1\in[B]^{<\omega}$ and such that $z_i\cap\bigcap\sigma_i\subseteq^*x_{0,\delta}$ for all $i<2$ and $\delta<\lambda$.  Hence, $\bigcap\bigcup_{i<2}\sigma_i\subseteq^*x_{0,\delta}$ for all $\delta<\lambda$, in contradiction with (\ref{eqlattice}).
\end{proof}

If $\alpha>0$, then choose $z$ as in the above claim; otherwise, choose $z$ arbitrarily.  If $\alpha=0$, then set $\rho_{\alpha,\beta}=\beta+1$.  Otherwise, we may choose $\rho_{\alpha,\beta}\in S_\alpha$ such that $\rho_{\alpha,\beta}>\rho_{\alpha,\gamma}$ for all $\gamma<\beta$ and
\begin{equation*}
\forces_{\lambda\alpha+\beta}\forall\sigma\in[A_{\rho_{\alpha,\beta}}]^{<\omega}\ \exists\delta<\rho_{\alpha,\beta}\ \,z\cap\bigcap\sigma\not\subseteq^*x_{0,\delta}.
\end{equation*} 
Set $D_{\alpha,\beta+1}=D_{\alpha,\beta}\cup\{\la\alpha,\rho_{\alpha,\beta}\ra\}$.  Let $A'$ be a $\mbP_{\lambda\alpha+\beta}$\nbd-name forced to satisfy $A'=A_{\rho_{\alpha,\beta}}\cup\{z\}$ if $z$ splits $B$ and $A'=A_{\rho_{\alpha,\beta}}$ otherwise.  Let $\mbQ_0$ be a name for the Booth forcing for $A'\cup\{\omega\setminus n:n<\omega\}$; let $x_{\alpha,\beta}$ be a name for a generic pseudointersection of $A'\cup\{\omega\setminus n:n<\omega\}$.  (The purpose of $\{\omega\setminus n:n<\omega\}$ is to ensure that $x_{\alpha,\beta}$ does not almost contain any element of $[\omega]^\omega\cap V^{\mbP_{\lambda\alpha+\beta}}$.)

Let $F_{\lambda\alpha+\beta}$ to be a $\mbP_{\lambda\alpha+\beta}$\nbd-name for a surjection from $\lambda$ to the elements of $[[\omega]^\omega]^{<\kappa}$ that have the SFIP.  \Wma\ that corresponding $F_\gamma$ have already been constructed for all $\gamma<\lambda\alpha+\beta$.  Let $\mbQ_1$ be a name for the Booth forcing for $F_{\lambda\alpha+\zeta}(\eta)$.

Further suppose we have constructed sequences $\la w_\gamma\ra_{\gamma<\alpha}$ and $\la U_\gamma\ra_{\gamma<\alpha}$ of $\mbP_{\lambda\alpha}$\nbd-names such that $\forces_{\lambda\gamma}U_\delta\cup\{w_\delta\}\subseteq U_\gamma\in\omega^*$ for all $\delta<\gamma<\alpha$, and such that $w_\gamma$ is forced to be a pseudointersection $U_\gamma$ for all $\gamma<\alpha$.  If $\beta\not=0$, then let $\mbQ_2$ be a name for the trivial forcing.  If $\beta=0$, then choose $U_\alpha$ such that $\forces_{\lambda\alpha}\bigcup_{\gamma<\alpha}U_\gamma\cup\{w_\gamma\}\subseteq U_\alpha\in\omega^*$, let $\mbQ_2$ be a name for the Booth forcing for $U_\alpha$, and let $w_\alpha$ be a name for a generic pseudointersection of $U_\alpha$.  

Further suppose we have constructed a sequence $\la a_\gamma\ra_{\gamma<\alpha}$ of $\mbP_{\lambda\alpha}$\nbd-names whose range is forced to be an almost disjoint subfamily of $[\omega]^\omega$.  If $\beta\not=0$, then let $\mbQ_3$ be a name for the trivial forcing.  If $\beta=0$, then let $\mbQ_3$ be a name for the Booth forcing for $\{\omega\setminus a_\gamma:\gamma<\alpha\}$, and let $a_\alpha$ be a name for a generic pseudointersection of $\{\omega\setminus a_\gamma:\gamma<\alpha\}$. 

Further suppose we have constructed a sequence $\la i_\gamma\ra_{\gamma<\alpha}$ of $\mbP_{\lambda\alpha}$\nbd-names whose range is forced to be an independent subfamily of $[\omega]^\omega$.  If $\beta\not=0$, then let $\mbQ_4$ be a name for the trivial forcing.  If $\beta=0$, then set $I=\{i_\gamma:\gamma<\alpha\}$ and let $J$ and $x$ be as in Lemma~\ref{LEMmif}; let $\mbQ_4$ be a name for the Booth forcing for $J$; let $i_\alpha$ be a name for $x$. 

Set $\mbP_{\lambda\alpha+\beta+1}=\mbP_{\lambda\alpha+\beta}*\prod_{n<5}\mbQ_n$.  \Wma\ $\card{\prod_{n<5}\mbQ_n}\leq\lambda$; hence, $\mbP_{\lambda\alpha+\beta+1}$ has property (K) and size at most $\lambda$.  Also, $B\cup\{x_{\alpha,\beta}\}$ is forced to have the SFIP by $\mbQ_0$\nbd-genericity because for every $b\in B$ we have that $\{b\}\cup A'$ is forced to have the SFIP because $\{b\}\cup A'\subseteq B\cup\{z\}$ if $z$ splits $B$ and $\{b\}\cup A'\subseteq B$ otherwise.  Let us also show that (\ref{eqlattice}) holds if we replace $\beta$ with $\beta+1$.  \Wma\ $\alpha>0$.  Let $\sigma\in[B]^{<\omega}$.  Then there exists $\delta<\lambda$ such that $\forces_{\lambda\alpha+\beta}z\cap\bigcap(\sigma\cup\tau)\not\subseteq^*x_{0,\delta}$ for all $\tau\in[A_{\rho_{\alpha,\beta}}]^{<\omega}$; hence, $\bigl\{\bigl(\bigcap\sigma\bigr)\setminus x_{0,\delta}\bigr\}\cup A'$ is forced to have the SFIP; hence, $\forces_{\lambda\alpha+\beta+1}x_{\alpha,\beta}\cap\bigcap\sigma\not\subseteq^*x_{0,\delta}$ by $\mbQ_0$\nbd-genericity.  Thus, (\ref{eqlattice}) holds as desired.

To complete our inductive construction of $\la\mbP_\gamma\ra_{\gamma\leq\lambda\kappa}$, it suffices to show that $\{\la\la\gamma,\rho_{\gamma,\delta}\ra,x_{\gamma,\delta}\ra:\la\gamma,\delta\ra\in E_{\alpha,\beta+1}\}$ is forced to be an order embedding of $D_{\alpha,\beta+1}$ into $\la[\omega]^\omega,\supseteq^*\ra$.  Suppose $\la\gamma,\delta\ra\in E_{\alpha,\beta}$.  Then $\la\alpha,\rho_{\alpha,\beta}\ra\not\leq\la\gamma,\rho_{\gamma,\delta}\ra$ and $\forces_{\lambda\alpha+\beta+1}x_{\alpha,\beta}\not\supseteq^*x_{\gamma,\delta}$ by $\mbQ_0$\nbd-genericity.  If $\la\gamma,\rho_{\gamma,\delta}\ra<\la\alpha,\rho_{\alpha,\beta}\ra$, then $x_{\gamma,\delta}\in A'$; whence, $\forces_{\lambda\alpha+\beta+1}x_{\gamma,\delta}\supsetneq^*x_{\alpha,\beta}$.  Suppose $\la\gamma,\rho_{\gamma,\delta}\ra\not<\la\alpha,\rho_{\alpha,\beta}\ra$.  Then $\rho_{\alpha,\beta}<\rho_{\gamma,\delta}$; hence, $\rho_{\gamma,\delta}\geq\rho_{\alpha,\beta}+1=\rho_{0,\rho_{\alpha,\beta}}$; hence, $x_{\gamma,\delta}\subseteq^*x_{0,\rho_{\alpha,\beta}}$.  By construction, $A'\cup\{\omega\setminus x_{0,\rho_{\alpha,\beta}}\}$ is forced to have the SFIP; hence, $\forces_{\lambda\alpha+\beta+1}x_{\gamma,\delta}\subseteq^*x_{0,\rho_{\alpha,\beta}}\not\supseteq^*x_{\alpha,\beta}$ by $\mbQ_0$\nbd-genericity.  Thus, $\{\la\la\gamma,\rho_{\gamma,\delta}\ra,x_{\gamma,\delta}\ra:\la\gamma,\delta\ra\in E_{\alpha,\beta+1}\}$ is forced to be an embedding as desired.

Let us show that $V^{\mbP_{\lambda\kappa}}$ satisfies $\lambda\leq\ochar{\omega^*}$. Let $G$ be a generic filter of $\mbP_{\lambda\kappa}$ and set $\mcB=\{(x_{\alpha,\beta})_G^*:\la\alpha,\beta\ra\in\kappa\times\lambda\}$.  Then $\mcB$ is a local base at some $p\in(\omega^*)^{V[G]}$ because every element of $([\omega]^\omega)^{V[G]}$ is handled by an appropriate $\mbQ_0$.  By Lemma~\ref{LEMmutuallydense}, $\mcB$ contains a \cardop{\ochar{p,\omega^*}}\ local base $\{(x_{\alpha,\beta})_G^*:\la\alpha,\beta\ra\in I\}$ at $p$ for some $I\subseteq\kappa\times\lambda$.  Set $J=\{\la\alpha,\rho_{\alpha,\beta}\ra:\la\alpha,\beta\ra\in I\}$.  Then $J$ is cofinal in $\kappa\times\lambda$; hence, by Lemma~\ref{LEMkappacrosslambda}, $J$ is not $\nu$\nbd-like for any $\nu<\lambda$.  Hence, $\ochar{\omega^*}^{V[G]}\geq\lambda$.

Finally, let us show that $V^{\mbP_{\lambda\kappa}}$ satisfies $\max\{\carda,\cardi,\cardu\}\leq\kappa\leq\cardp$.   Working in $V[G]$, notice that $\cardu\leq\kappa$ because $\bigcup_{\alpha<\kappa}(U_\alpha)_G\in\omega^*$ and $\{(w_\alpha)_G^*:\alpha<\kappa\}$ is a local base at $\bigcup_{\alpha<\kappa}(U_\alpha)_G$.  Moreover, $\{(a_\alpha)_G:\alpha<\kappa\}$ and $\{(i_\alpha)_G:\alpha<\kappa\}$ witness that $\carda\leq\kappa$ and $\cardi\leq\kappa$.  For $\cardp\geq\kappa$, note that very element of $[[\omega]^\omega]^{<\kappa}$ with the SFIP is $\left(F_{\lambda\alpha+\zeta}(\eta)\right)_G$ for some $\alpha<\kappa$ and $\zeta,\eta<\lambda$. By $\mbQ_1$\nbd-genericity, a pseudointersection of $\left(F_{\lambda\alpha+\zeta}(\eta)\right)_G$ is added at stage $\lambda\alpha+\varphi(\zeta,\eta)$.
\end{proof}

\begin{theorem}
$\opichar{\omega^*}=\omega$.
\end{theorem}
\begin{proof}
Fix $p\in\omega^*$.  By a result of Balcar and Vojt\'a\v s~\cite{balcarvojtas}, there exists $\la y_x\ra_{x\in p}$ such that $y_x\in[x]^{\omega}$ for all $x\in p$ and $\{y_x\}_{x\in p}$ is an almost disjoint family.  Clearly, $\{y_x^*\}_{x\in p}$ is a pairwise disjoint---and therefore \omegaop---local $\pi$\nbd-base at $p$.
\end{proof}

\section{Powers of $\omega^*$}

\begin{definition}
A \emph{box} is a subset $E$ of a product space $\prod_{i\in I}X_i$ such that there exist $\sigma\in[I]^{<\omega}$ and $\la E_i\ra_{i\in\sigma}$ such that $E=\bigcap_{i\in\sigma}\inv{\pi_i}E_i$.  Let $\owbox{\prod_{i\in I}X_i}$ denote the least infinite $\kappa$ such that $\prod_{i\in I}X_i$ has a \kappaop\ base of open boxes.
\end{definition}

\begin{lemma}[Peregudov~\cite{peregudov97}]\label{LEMowprodsup}
In any product space $X=\prod_{i\in I}X_i$, we have $\ow{X}\leq\owbox{X}\leq\sup_{i\in I}\ow{X_i}$.
\end{lemma}

\begin{lemma}[Malykhin~\cite{malykhin}]\label{LEMowbigprod}
Let $X=\prod_{i\in I}X_i$ where each $X_i$ is a nonsingleton $T_1$ space.  If $\weight{X}\leq\card{I}$, then $\ow{X}=\owbox{X}=\omega$.
\end{lemma}

\begin{remark}
In Lemma~\ref{LEMowbigprod}, the hypothesis that the factor spaces be nonsingleton and $T_1$ can be weakened to merely require that each factor space is the union of two nontrivial open sets.  Also, the conclusion of Lemma~\ref{LEMowbigprod} may be amended with the statement that $X$ has a $\la\card{I},\omega\ra$\nbd-splitter: use $\la\{\inv{\pi}_i U_i,\inv{\pi}_i V_i\}\ra_{i\in I}$ where each $\{U_i,V_i\}$ is a nontrivial open cover of $X_i$.
\end{remark}

\begin{theorem}\label{THMnonincreaseomega}
The sequence $\la\ow{(\omega^*)^{\omega+\alpha}}\ra_{\alpha\in\ord}$ is nonincreasing.  Moreover, $\ow{(\omega^*)^\cardc}=\omega$.
\end{theorem}
\begin{proof}
Note that if $\omega\leq\alpha\leq\beta$, then $(\omega^*)^\beta\homeo((\omega^*)^\alpha)^\beta$.  Then apply Lemmas~\ref{LEMowprodsup} and \ref{LEMowbigprod}.
\end{proof}

\begin{lemma}\label{LEMowbox}
Let $0<n<\omega$ and $X$ be a space.  Then $\owbox{X^n}=\ow{X}$.
\end{lemma}
\begin{proof}
Set $\kappa=\owbox{X^n}$.  By Lemma~\ref{LEMowprodsup}, $\kappa\leq\ow{X}$.  Let us show that $\ow{X}\leq\kappa$.  Let $\mcA$ be a \kappaop\ base of $X^n$ consisting only of boxes.  Let $\mcB$ denote the set of all nonempty open $V\subseteq X$ for which there exists $\prod_{i<n}U_i\in\mcA$ such that $V=\bigcap_{i<n}U_i$.  Then $\mcB$ is a base of $X$ because if $p\in U$ and $U$ is an open subset of $X$, then there exists $\prod_{i<n}U_i\in\mcA$ such that $\la p\ra_{i<n}\in\prod_{i<n}U_i\subseteq U^n$; whence, $p\in\bigcap_{i<n}U_i\subseteq U$ and $\bigcap_{i<n}U_i\in\mcB$.  

It suffices to show that $\mcB$ is \kappaop.  Suppose not.  Then there exist $\prod_{i<n}U_i\in\mcA$ and $\la\prod_{i<n}V_{\alpha,i}\ra_{\alpha<\kappa}\in\mcA^\kappa$ such that 
\begin{equation*}
\emptyset\not=\bigcap_{i<n}U_i\subseteq\bigcap_{i<n}V_{\alpha,i}\not=\bigcap_{i<n}V_{\beta,i}
\end{equation*}
 for all $\alpha<\beta<\kappa$.  Clearly, $\prod_{i<n}V_{\alpha,i}\not=\prod_{i<n}V_{\beta,i}$ for all $\alpha<\beta<\kappa$.  Choose $U\in\mcA$ such that $U\subseteq(\bigcap_{i<n}U_i)^n$.  Then $U\subseteq\prod_{i<n}V_{\alpha,i}$ for all $\alpha<\kappa$, in contradiction with how we chose $\mcA$.
\end{proof}

\begin{lemma}\label{LEMowfinprodweightchar}
If $0<n<\omega$ and $X$ is a compact space such that $\character{p,X}=\weight{X}$ for all $p\in X$, then $\ow{X}=\ow{X^n}$.
\end{lemma}
\begin{proof}
By Lemma~\ref{LEMowbox}, it suffices to show that $\owbox{X^n}\leq\ow{X^n}$.  By Lemma~\ref{LEMowlowerboundcharweight}, either $X^n$ has a $\la\weight{X^n},\ow{X^n}\ra$\nbd-splitter, or $\ow{X^n}=\weight{X^n}^+$.  Hence, by Lemma~\ref{LEMowsplitterupperbound}, $\owbox{X^n}\leq\ow{X^n}$.
\end{proof}

\begin{theorem}\label{THMowomegastarfinprod}
If $0<n<\omega$, then $\ow{\omega^*}\geq\ow{(\omega^*)^n}\geq\min\{\ow{\omega^*},\cardc\}$.  Moreover, $\max\{\cardu,\cf\cardc\}=\cardc$ implies $\ow{\omega^*}=\ow{(\omega^*)^n}$.
\end{theorem}
\begin{proof}
Lemma~\ref{LEMowprodsup} implies $\ow{\omega^*}\geq\ow{(\omega^*)^n}$.  To prove the rest of the theorem, first consider the case $\cardr<\cardc$.  As in the proof of Theorem~\ref{THMowomegastarcplus}, construct a point $p\in\omega^*$ such that $\pi\character{p,\omega^*}=\cardr$ and $\character{p,\omega^*}=\cardc$.  Then $\pi\character{\la p\ra_{i<n},(\omega^*)^n}=\cardr$ and $\character{\la p\ra_{i<n},(\omega^*)^n}=\cardc$; hence, $\ow{(\omega^*)^n}\geq\cardc$ by Theorem~\ref{PROpicharcfchar}.  Moreover, if $\cf\cardc=\cardc$, then $\ow{(\omega^*)^n}=\ow{\omega^*}=\cardc^+$.  If $\cardu=\cardc$, then $\ow{\omega^*}=\ow{(\omega^*)^n}$ by Lemma~\ref{LEMowfinprodweightchar}.  Finally, in the case $\cardr=\cardc$, we have $\cardu=\cardc$, which again implies $\ow{\omega^*}=\ow{(\omega^*)^n}$.
\end{proof}

\begin{corollary}\label{CORnonincreasing}
Suppose $\max\{\cardu,\cf\cardc\}=\cardc$.  Then $\la\ow{(\omega^*)^{1+\alpha}}\ra_{\alpha\in\ord}$ is nonincreasing.
\end{corollary}
\begin{proof}
By Theorem~\ref{THMowomegastarfinprod} and Lemma~\ref{LEMowprodsup}, $\ow{(\omega^*)^n}=\ow{\omega^*}\geq\ow{(\omega^*)^\alpha}$ whenever $0<n<\omega\leq\alpha$.  The rest follows from Theorem~\ref{THMnonincreaseomega}.
\end{proof}

\begin{theorem}
Suppose $\cardu=\cardc$.  Then $\ow{(\omega^*)^{1+\alpha}}=\ow{\omega^*}$ for all $\alpha<\cf\cardc$.
\end{theorem}
\begin{proof}
Let $\lambda$ be an arbitrary infinite cardinal less than $\ow{\omega^*}$.  By Lemma~\ref{LEMowlowerboundcharweight}, it suffices to show that $(\omega^*)^{1+\alpha}$ does not have a $\la\cardc,\lambda\ra$\nbd-splitter.  Seeking a contradiction, suppose $\la\mcF_\beta\ra_{\beta<\cardc}$ is such a $\la\cardc,\lambda\ra$\nbd-splitter.  \Wma\ $\bigcup_{\beta<\cardc}\mathcal{F}_\beta$ consists only of open boxes because we can replace each $\mcF_\beta$ with a suitable refinement.  Since $\alpha<\cf\cardc$, there exist $\sigma\in[1+\alpha]^{<\omega}$ and $I\in[\cardc]^\cardc$ such that, for every $U\in\bigcup_{\beta\in I}\mcF_\beta$, there exists $\varphi(U)\subseteq(\omega^*)^\sigma$ such that $U=\inv{\pi}_\sigma\varphi(U)$.  Let $j$ be a bijection from $\cardc$ to $I$.  Then $\la\varphi``\mcF_{j(\beta)}\ra_{\beta<\cardc}$ is a $\la\cardc,\lambda\ra$\nbd-splitter of $(\omega^*)^\sigma$.  Hence, $\ow{(\omega^*)^\sigma}\leq\lambda<\ow{\omega^*}$ by Lemma~\ref{LEMowsplitterupperbound}.  But $\ow{(\omega^*)^\sigma}<\ow{\omega^*}$ contradicts Theorem~\ref{THMowomegastarfinprod}.
\end{proof}

\begin{lemma}\label{LEMcfwcfw}
Suppose a space $X$ has a $\la\cf\weight{X},\cf\weight{X}\ra$\nbd-splitter.  Then $\ow{X}\leq\weight{X}$.
\end{lemma}
\begin{proof}
Set $\kappa=\cf\weight{X}$ and $\lambda=\weight{X}$.  Let $\la\mcF_\alpha\ra_{\alpha<\kappa}$ be a $\la\kappa,\kappa\ra$\nbd-splitter of $X$.  Let $h:\lambda\rightarrow\kappa$ satisfy $\card{\inv{h}\{\alpha\}}<\lambda$ for all $\alpha<\kappa$.  Then $\la\mcF_{h(\alpha)}\ra_{\alpha<\lambda}$ is a $\la\lambda,\lambda\ra$\nbd-splitter because if $I\in[\lambda]^\lambda$, then $h``I\in[\kappa]^\kappa$.  By Lemma~\ref{LEMowsplitterupperbound}, $\ow{X}\leq\lambda$.
\end{proof}

\begin{remark}
The proof of the above lemma shows that for any infinite cardinal $\kappa$, a space with a $\la\cf\kappa,\cf\kappa\ra$\nbd-splitter also has a $\la\kappa,\kappa\ra$\nbd-splitter.
\end{remark}

\begin{theorem}\label{THMowprodcf}
$\ow{(\omega^*)^{\cf\cardc}}\leq\cardc$.
\end{theorem}
\begin{proof}
The sequence $\la\{\inv{\pi}_\alpha(\{2n:n<\omega\}^*),\inv{\pi}_\alpha(\{2n+1:n<\omega\}^*)\}\ra_{\alpha<\cf\cardc}$ is a $\la\cf\cardc,\omega\ra$\nbd-splitter of $(\omega^*)^{\cf\cardc}$.  Apply Lemma~\ref{LEMcfwcfw}.
\end{proof}

\begin{theorem}
For all cardinals $\kappa$ satisfying $\kappa>\cf\kappa>\omega_1$, it is consistent that $\cardc=\kappa$ and $\cardr<\cf\cardc$.  The last inequality implies $\ow{(\omega^*)^{1+\alpha}}=\cardc^+$ for all $\alpha<\cf\cardc$ and $\ow{(\omega^*)^\beta}=\cardc=\kappa$ for all $\beta\in\cardc\setminus\cf\cardc$.
\end{theorem}
\begin{proof}
Starting with $\cardc=\kappa$ in the ground model, the proof of Theorem~\ref{THMowomegastarcplus} shows how to force $\cardr=\cardu=\omega_1$ while preserving $\cardc$.  Now suppose $\cardr<\cf\cardc$.  Fix $\alpha<\cf\cardc$ and $\beta\in\cardc\setminus\cf\cardc$.  By Theorems~\ref{THMowprodcf} and \ref{THMnonincreaseomega}, $\ow{(\omega^*)^\beta}\leq\cardc$.  To see that $\ow{(\omega^*)^\beta}\geq\cardc$, proceed as in the proof of Theorem~\ref{THMowomegastarfinprod}, constructing a point with character $\cardc$ and $\pi$\nbd-character $\card{\beta}$.  Similarly prove $\ow{(\omega^*)^{1+\alpha}}=\cardc^+$ by constructing a point with character $\cardc$ and $\pi$\nbd-character $\card{\cardr+\alpha}$.
\end{proof}

\begin{lemma}\label{LEMowprodgmusets}
Suppose $\kappa$, $\lambda$, and $\mu$ are cardinals and $p$ is a point in a product space $X=\prod_{\alpha<\kappa}X_\alpha$ satisfying the following for all $\alpha<\kappa$.
\begin{enumerate}
\item\label{enumkappalambdasmall1} $0<\kappa<\weight{X}$ and $\omega\leq\lambda\leq\weight{X}$.
\item\label{enumkappalambdasmall2} $\kappa<\cf\weight{X}$ or $\lambda<\weight{X}$.
\item\label{enumlambdamusmall} $\mu<\lambda$ or $\mu=\cf\lambda$.
\item\label{enumGmuset} $\character{p(\alpha),X_\alpha}<\lambda$ or the intersection of any $\mu$\nbd-many neighborhoods of $p(\alpha)$ has nonempty interior.
\end{enumerate}
Then $\character{p,X}<\weight{X}$ or $\ow{X}>\mu$.
\end{lemma}
\begin{proof}
Let $\mcA$ be a base of $X$.  Set $\mcB=\{U\in\mcA:p\in U\}$.  For each $\alpha<\kappa$, let $\mcC_\alpha$ be a local base at $p(\alpha)$ of size $\character{p(\alpha),X_\alpha}$.  Set $F=\bigcup_{r\in[\kappa]^{<\omega}}\prod_{\alpha\in r}\mcC_\alpha$.  For each $\sigma\in F$, set $U_\sigma=\bigcap_{\alpha\in\dom\sigma}\inv{\pi_\alpha}\sigma(\alpha)$.  For each $V\in\mcB$, choose $\sigma(V)\in F$ such that $p\in U_{\sigma(V)}\subseteq V$. \Wma\ $\character{p,X}=\weight{X}$; hence, by (\ref{enumkappalambdasmall1}) and (\ref{enumkappalambdasmall2}), there exist $r\in[\kappa]^{<\omega}$ and $\mcD\in[\mcB]^\lambda$ such that $\dom\sigma(V)=r$ for all $V\in\mcD$.  Set $s=\{\alpha\in r:\character{p(\alpha),X_\alpha}<\lambda\}$ and $t=r\setminus s$.  By (\ref{enumlambdamusmall}), there exist $\tau\in\prod_{\alpha\in s}\mcC_\alpha$ and $\mcE\in[\mcD]^{\mu}$ such that $\sigma(V)\restrict s=\tau$ for all $V\in\mcE$.  By (\ref{enumGmuset}), $\bigcap_{V\in\mcE}\sigma(V)(\alpha)$ has nonempty interior for all $\alpha\in t$.  Hence, $\bigcap\mcE$ has nonempty interior because it contains $U_\tau\cap\bigcap_{\alpha\in t}\inv{\pi_\alpha}\bigcap_{V\in\mcE}\sigma(V)(\alpha)$.  Thus, $\ow{X}>\mu$.
\end{proof}

\begin{theorem}\label{THMowsmallprodcardp}
Suppose $0<\alpha<\cardc$ and $\la X_\beta\ra_{\beta<\alpha}$ is a sequence of spaces each with weight at most $\cardc$.  Then $\ow{\prod_{\beta<\alpha}(X_\beta\oplus\omega^*)}>\nu$ for all regular $\nu<\cardp$.
\end{theorem}
\begin{proof}
Let $\nu$ be an arbitrary infinite regular cardinal less than $\cardp$. Set $\kappa=\card{\alpha}$ and $\lambda=\mu=\nu$.  Choose $q\in\omega^*$ such that $\character{q,\omega^*}=\cardc$; set $p=\la q\ra_{\beta<\alpha}$.  Applying Lemma~\ref{LEMowprodgmusets}, we have $\ow{\prod_{\beta<\alpha}(X_\beta\oplus\omega^*)}>\nu$.
\end{proof}

\begin{corollary}
Suppose $\cardp=\cardc$.  Then $\ow{(\omega^*)^{1+\alpha}}=\cardc$ for all $\alpha<\cardc$.
\end{corollary}
\begin{proof}
By Theorem~\ref{THMowcardr}, $\ow{\omega^*}\leq\cardc$.  Hence, by Corollary~\ref{CORnonincreasing}, $\ow{(\omega^*)^{1+\alpha}}\leq\cardc$ for all $\alpha\in\ord$.  By Theorem~\ref{THMowsmallprodcardp}, $\ow{(\omega^*)^{1+\alpha}}\geq\cardc$ for all $\alpha<\cardc$.
\end{proof}

\begin{corollary}
Suppose $\alpha<\cardc$ and $\la X_\beta\ra_{\beta<\alpha}$ is a sequence of spaces each with weight at most $\cardc$.  Then $\prod_{\beta<\alpha}(X_\beta\oplus\omega^*)$ is not homeomorphic to a product of $\cardc$\nbd-many nonsingleton spaces.
\end{corollary}
\begin{proof}
Combine Theorem~\ref{THMowsmallprodcardp} and Lemma~\ref{LEMowbigprod}.
\end{proof}

\section{Questions}

\begin{question}
Is it consistent that $\ow{\omega^*}=\cardc^+$ and $\cardr\geq\cf\cardc$?
\end{question}

\begin{question}
Is $\ow{\omega^*}<\cardss_\omega$ consistent?  This inequality implies $\cardu<\cardc$.  Hence, by Theorem~\ref{THMowcardr}, the inequality further implies 
\begin{equation*}
\cf\cardc\leq\cardr\leq\cardu<\cardc=\ow{\omega^*}<\cardss_\omega=\cardc^+.
\end{equation*}
More generally, does any space $X$ have a base that does not contain an \cardop{\ow{X}} base?
\end{question}

\begin{question}
Is $\cardss_\omega<\cardss_2$ consistent?
\end{question}

\begin{question}
Letting $\cardg$ denote the groupwise density number (see 6.26 of \cite{blass}), is $\ow{\omega^*}<\cardg$ consistent? $\ochar{\omega^*}<\cardg$?  In particular, what are $\ow{\omega^*}$ and $\ochar{\omega^*}$ in the Laver model (see 11.7 of \cite{blass})?
\end{question}

\begin{question}
Is $\cf\ow{\omega^*}<\ow{\omega^*}<\cardc$ consistent? $\cf\ow{\omega^*}=\omega$?
\end{question}

\begin{question}
Is $\cf\cardc<\ow{\omega^*}<\cardc$ consistent?
\end{question}

\begin{question}
What is $\ochar{\omega^*}$ in the forcing extension of the proof of Theorem~\ref{THMcardblessowomegastar}?  More generally, is it consistent that $\ochar{\omega^*}<\ow{\omega^*}\leq\cardc$?
\end{question}

\begin{question}
Is $\ochar{\omega^*}=\omega$ consistent?  An affirmative answer would be a strengthening of Shelah's result~\cite{shelah} that $\omega^*$ consistently has no P\nbd-points.  If the answer is negative, then which, if any, of $\cardp$, $\cardh$, $\cards$, and $\cardg$ are lower bounds of $\ochar{\omega^*}$ in ZFC?
\end{question}

\begin{question}
Is $\cf\cardc<\ochar{\omega^*}$ consistent?  $\cf\cardc<\ochar{\omega^*}<\cardc$?
\end{question}

\begin{question}
Does any Hausdorff space have uncountable local Noetherian $\pi$\nbd-type?  (It is easy to construct such $T_1$ spaces: give $\omega_1+1$ the topology $\{(\omega_1+1)\setminus(\alpha\cup\sigma):\alpha<\omega_1\text{ and }\sigma\in[\omega_1+1]^{<\omega}\}\cup\{\emptyset\}$.)
\end{question}

\begin{question}
Is it consistent that $\ow{(\omega^*)^{1+\alpha}}<\min\{\ow{\omega^*},\cardc\}$ for some $\alpha<\cardc$?
Is it consistent that $\ow{(\omega^*)^{1+\alpha}}<\ow{\omega^*}$ for some $\alpha<\cf\cardc$?
\end{question}

\end{document}